\documentclass[12pt,a4]{article}
\usepackage{natbib}
\usepackage{amsthm,amsfonts,epsfig,color}
\usepackage{mathptmx}
\usepackage{bm}
\usepackage{amsmath}

\usepackage{color}
\usepackage[normalem]{ulem}
\usepackage[T1]{fontenc}
\usepackage{epsfig}
\usepackage{graphicx}
\usepackage{graphics}
\usepackage{pdflscape}
\usepackage{amssymb,srcltx,amsthm,epsfig}
\usepackage{appendix}
\newcommand{\lp}{\left(}
\newcommand{\rp}{\right)}

\newcommand{\sumi}{\sum^n_{i=1}}

\newcommand{\sumas}{\sum^n_{i=1}}

\newcommand{\is}{i=1,\ldots,n}

\newcommand{\by}{{\bf y}}
\newcommand{\yp}{\mathbf{y}}
\newcommand{\Y}{\mathbf{Y}}

\newcommand{\X}{\mathbf{X}}

\newcommand{\y}{\mathbf{y}}

\newcommand{\bu}{\mathbf{u}}
\newcommand{\bx}{\mathbf{x}}
\newcommand{\xp}{\mathbf{x}}

\newcommand{\bbeta}{\mbox{${\bm \beta}$}}
\newcommand{\btheta}{\mbox{${ \bm \theta}$}}

\newcommand{\bxi}{\mbox{${ \bm \xi}$}}

\newtheorem{theorem}{Theorem}[section]
\begin{document}

\title{ Case-Deletion Diagnostics  for Quantile Regression Using the Asymmetric Laplace Distribution}

\author{\small  Luis E. Benites$^{a}$  \and \small Victor H. Lachos$^{b}$\thanks{Corresponding author.  Address for correspondence: Departamento de Estat\'{\i}stica,  Rua S\'ergio Buarque de Holanda,
651,  Cidade Universit\'aria Zeferino Vaz, Campinas, S\~ao Paulo,
Brazil. CEP 13083-859.  E-mail: \texttt{hlachos@ime.unicamp.br}} \and  \small Filidor E. Vilca$^{b}$ \\
{\em \small (a) Departamento de Estat\'{\i}stica, Universidade de S\~{a}o Paulo, Brazil }  \\
{\em \small (b) Departamento de Estat\'{\i}stica, Universidade Estadual de Campinas, Brazil }  \\
}

\date{}
\maketitle
\begin{abstract}
To make inferences about the shape of a population distribution, the
widely popular mean regression model, for example, is inadequate if
the distribution is not approximately Gaussian (or symmetric).
Compared to conventional mean regression (MR), quantile regression
(QR) can characterize the entire conditional distribution of the
outcome variable, and is more robust to outliers and
misspecification of the error distribution. We present a
likelihood-based approach to the estimation of the regression
quantiles based on the asymmetric Laplace distribution (ALD), which
has  a hierarchical representation that facilitates the
implementation of the EM algorithm for the  maximum-likelihood
estimation. We develop  a case-deletion diagnostic analysis for QR
models   based on the conditional expectation of the complete-data
log-likelihood function related to the EM algorithm.  The techniques
are illustrated with both simulated and real data sets, showing that
our approach out-performed other common classic estimators. The
proposed algorithm and methods are implemented in the R package
\verb"ALDqr()".
\vspace*{0.5cm}\\
\noindent {\bf Keywords}: Quantile regression model; EM algorithm;
Case-deletion model; Asymmetric Laplace distribution.
\end{abstract}
\section{Introduction}

QR models have become increasingly popular since the seminal work of
\cite{koenker1978regression}. In contrast to the mean regression
model, QR belongs to a robust model family, which can give an
overall assessment of the covariate effects at different quantiles
of the outcome \citep{koenker2005quantile}. In particular, we can
model the lower or higher quantiles of the outcome to provide a
natural assessment of covariate effects specific for those
regression quantiles. Unlike conventional models, which only address
the conditional mean or the central effects of the covariates, QR
models quantify the entire conditional distribution of the outcome
variable. In addition, QR does not impose any distributional
assumption on the error, except requiring  the error to have a zero
conditional quantile. The foundations of the methods for independent
data are now consolidated, and some statistical methods for
estimating and drawing inferences about conditional quantiles are
provided by most of the available statistical programs (e.g., R,
SAS, Matlab and Stata). For instance, just to name a few, in the
well-known R package \verb"quantreg()" is implemented  a variant of
the \cite{barrodale1977algorithms} simplex (BR) for linear
programming problems described in \cite{koenker1987algorithm}, where
the standard errors are computed by the rank inversion method
\citep{koenker2005quantile}. Another method implemented in this
popular package is Lasso Penalized Quantile Regression (LPQR),
introduced by \cite{tibshirani1996regression}, where a penalty
parameter is specified to determine how much shrinkage occurs in the
estimation process. QR can be implemented in a range of different
ways. \cite{koenker2005quantile} provided an overview of some
commonly used quantile regression techniques from a
"classical" framework.

\cite{kottas2001bayesian} considered
median regression from a Bayesian point of view, which is a special case of quantile regression,
and discussed non-parametric modeling for the error distribution
based on either P\'{o}lya tree or Dirichlet process priors. Regarding
general quantile regression, \cite{yu2001bayesian} proposed a
Bayesian modeling approach by using the ALD,
\cite{kottas2009bayesian} developed Bayesian semi-parametric models
for quantile regression using Dirichlet process mixtures for the
error distribution, {\cite{geraci07} studied quantile regression
for longitudinal data using the ALD.} Recently, \cite{kozumi2011gibbs} developed
a simple and efficient Gibbs sampling algorithm for fitting the quantile regression model based on a location-scale mixture representation of the ALD.

An interesting aspect to be considered in statistical modelling is
the diagnostic analysis.  This can be carried out by conducting an
influence analysis for detecting influential observations.  One   of
the   most technique to detect influential observations is  the
case-deletion approach. The famous approach of Cook (1977)  has been
applied extensively to assess the influence of an observation in
fitting a statistical model; see \cite{cook82} and the references
therein.   It is difficult to apply this approach directly to the QR
model  because the underlying observed-data likelihood function is
not differentiable at zero.  \cite{zhu2001case} presents an approach
to perform diagnostic analysis for general statistical models that
is  based on  the  Q-displacement function. This approach  has  been
applied successfully to perform influence analysis in several
regression models, for example, \cite{xie2007case}  considered  in
multivariate  $t$   distribution,
\cite{Matos.Lachos.Bala.Labra.2012} obtained case-deletion measures
for mixed-effects models following the \cite{zhu2001case}'s approach
and   in  \cite{Zeller.Labra.Lachos.Balakrishnan.2010}  we  can  see
some  results about   local influence for mixed-effects models
obtained by using  the Q-displacement function.

Taking advantage of the likelihood structure imposed by the ALD, the
hierarchical representation of the ALD,   we  develop here  an
EM-type algorithm for obtaining the ML estimates at the $p$th level,
and  by simulation studies our EM algorithm outperformed the
competing  BR and LPQR algorithms, where  the standard error is
obtained as a by-product. Moreover, we obtain    case-deletion
measures for the QR model. Since QR methods complement and improve
established means regression models, we feel that the assessment of
robustness aspects of the parameter estimates in QR is also an
important concern at a given quantile level $p\in(0, 1)$.

\indent The rest of the paper is organized as follows. Section 2
introduces the connection between QR and ALD as well as outlining
the main results related to ALD.  Section 3 presents  an EM-type
algorithm to proceed with ML estimation for the parameters at the
$p$th level. Moreover,  the observed information matrix  is derived.
Section \ref{Sec Diagnostic}  provides a brief sketch  of the case-deletion method
for the model with incomplete data, and also develop a methodology pertinent to the ALD. Sections \ref{sec application}
and \ref{sec simulation study} are dedicated to the analysis of real
and simulated data sets, respectively. Section 6 concludes with a
short discussion of issues raised by our study and some possible
directions for the future research.

\section{The quantile regression model} \label{sec tCR}
Even though considerable amount of work has been done on  regression  models  and their extensions,  regression  models   by   using asymmetric Laplace distribution    have received little attention in the literature. Only recently, the a   study on quantile regression  model   based  on asymmetric Laplace distribution
 was presented by  Tian et al. (2014) who a derived several interesting
and attractive properties and   presented   an EM  algorithm. Before
presenting our derivation, let us recall firstly  the definition of
the asymmetric Laplace distribution and   after  this,  we    will
present  the  quantile regression model.

\subsection{Asymmetric Laplace distribution}
As discussed in \cite{yu2001bayesian}, we say that a random variable
Y is distributed as an ALD with location parameter $\mu$, scale
parameter $\sigma>0$ and skewness parameter $p\in (0,1)$, if its
probability density function (pdf) is given by
\begin{equation}\label{pdfAL}
f(y|\mu,\sigma,p)=\frac{p(1-p)}{\sigma}\exp\Big
\{-\rho_p\big(\frac{y-\mu}{\sigma}\big)\Big\},
\end{equation}
where $\rho_p(.)$ is the so called check (or loss) function defined
by $\rho_p(u)=u(p-\mathbb{I}\{u<0\})$, with   $\mathbb{I}\{.\}$
denoting the usual indicator function.   This  distribution is
denoted by $ALD(\mu,\sigma,p)$. It is easy to see that
$W=\rho_p\big(\frac{Y-\mu}{\sigma}\big)$ follows an exponential
distribution $\exp(1)$.

\begin{figure}[!htb]
\begin{center}
{\includegraphics[scale=0.6]{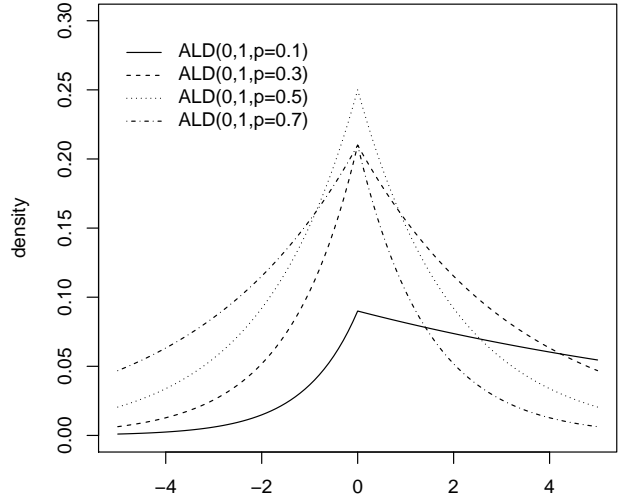}
\caption{Standard asymmetric Laplace density \label{fig:ald}}}
\end{center}
\end{figure}

A  stochastic representation   is  useful to  obtain  some  properties  of
the  distribution, as  for example, the moments, moment generating function (mgf), and estimation algorithm.  For the ALD     \cite{kotz2001laplace},  \cite{kozobowski00} and \cite{zhou13} presented  the  following  stochastic representation:
Let $U\sim
{\exp}(\sigma)$ and $Z\sim N(0,1)$ be two independent random variables. Then,
$Y\sim ALD(\mu,\sigma,p)$ can be   represented as
\begin{equation}\label{st-ALD}
Y\buildrel d\over=\mu+\vartheta_p U+\tau_p\sqrt{\sigma U} Z,
\end{equation}
where $\vartheta_p=\frac{1-2p}{p(1-p)}$ and
$\tau^2_p=\frac{2}{p(1-p)}$,  and $\buildrel d\over=$
denotes equality in distribution.
\indent Figure \ref{fig:ald}  shows how the skewness of the ALD
changes with altering values for $p$. For example, for $p=0.1$
almost all the mass of the ALD is situated in the right tail. For $p=0.5$, both tails of the ALD have equal mass and the  distribution then equals the more common double exponential
distribution. In contrast to the normal distribution with a
quadratic term in the exponent, the ALD is linear in the exponent.
This results in a more peaked mode for the ALD together with thicker
tails. On the other hand, the normal distribution has heavier
shoulders compared to the ALD.
From (\ref{st-ALD}),  we   have the hierarchical representation of   the ALD,  see  \cite{lum12},    given  by
\begin{eqnarray}
Y|U=u&\sim& N(\mu+\vartheta_p u,\tau^2_p\sigma u ),\label{hierar1}\\
U&\sim& exp(\sigma)\label{hierar2}.
\end{eqnarray}
This  representation   will   be   useful   for the  implementation  of  the  EM  algorithm. Moreover,  since   $Y|U=u\sim N(\mu+\vartheta_p u,\tau^2_p\sigma u )$, then one  can  derive   easily   the  pdf  of $Y$. That is,  the   pdf in  ( \ref{pdfAL})   can  be   expressed   as
\begin{equation}\label{pdfALs}
f(y|\mu,\sigma,p)=\frac{1}{\sqrt{2\pi}}\frac{1}{\tau_p\sigma^{\frac{3}{2}}} \exp\Big(\frac{\delta(y)}{\gamma}\Big)
A(y),
\end{equation}
where $\delta(y)=\frac{|y-\mu|}{\tau_p\sqrt{\sigma}}$, $\gamma=\sqrt{\frac{1}{\sigma}\big(2+\frac{\vartheta_p^2}{\tau^2_p}\big)}=\frac{\tau_p}{2\sqrt{\sigma}}$ and   $A(y)=2\Big(\frac{\delta(y)}{\gamma}\Big)^{1/2} K_{1/2}\big(\delta(y)\gamma\big)$, with $K_{\nu}(.)$  being  the  modified Bessel function of the third kind.  It   easy  to  see  that
that the conditional distribution of $U$, given $Y=y$, is
$U|(Y=y)\sim GIG (\frac{1}{2},\delta,\gamma)$. Here,  $GIG(\nu, a, b)$ denotes the Generalized Inverse
Gaussian (GIG) distribution; see \cite{barndorff2001non} for more details. The   pdf of GIG   distribution  is   given by
$$h(u|\nu,a,b)=\frac{(b/a)^{\nu}}{2K_{\nu}(ab)}u^{\nu-1}\exp\Big\{-\frac{1}{2}\big(a^2/{u}+b^2u\big)\Big\},\,\,u>0,\,\,\,\,\nu\in \mathbb{R},\,\,a,b>0.$$
The moments of $U$  can  be  expressed  as
$$E[U^k]=\left(\frac{a}{b}\right)^{k}\frac{K_{\nu+k}(ab)}{K_{\nu}(ab)},\,\,\ k\in \mathbb{R}.$$
Some   properties of the Bessel function of  the third kind
$K_{\lambda}(u)$    that will  be  useful  for the   developments
here  are: (i) $K_{\nu}(u)=K_{-\nu}(u)$; (ii)
$K_{\nu+1}(u)=\frac{2\nu}{u}K_{\nu}(u)+K_{\nu-1}(u)$; (iii)  for
non-negative integer  $r$,
$K_{r+1/2}(u)=\sqrt{\frac{\pi}{2u}}\exp(-u)\sum_{k=0}^r
\frac{(r+k)!(2u)^{-k}}{(r-k)!k!}$. A  special   case  is  $K_{1/2}(u)=
\sqrt{\frac{\pi}{2u}}\exp(-u)$.\\

\subsection{ Linear quantile regression}

Let $y_i$ be a response variable and $\bx_i$ a $k\times 1$
vector of covariates for the $i$th observation, and let
$Q_{y_i}(p|\bx_i)$ be the $p$th $(0 < p < 1)$ quantile regression
function of $y_i$ given $\bx_i$, $\is$ . Suppose that the relationship
between $Q_{y_i}(p|\bx_i)$ and $\bx_i$ can be modeled as $Q_{y_i}(p|\bx_i)
= \bx^{\top}_i\bbeta_p$, where $\bbeta_p$ is a vector $(k\times 1)$
of unknown parameters of interest. Then, we consider the quantile
regression model given by
\begin{equation}\label{QRmodel}
y_i = \bx^{\top}_i\bbeta_p + \epsilon_i,\,\,\,\is,
\end{equation}
where $\epsilon_i$ is the error term whose distribution (with density,
say, $f_p(.)$) is restricted to have the $p$th quantile equal to
zero, that is, $\int^{0}_{-\infty}f_p(\epsilon_i)d\epsilon_i=p$.  The error density $f_p(.)$ is often left unspecified in the classical literature. Thus, quantile regression estimation for
$\bbeta_p$ proceeds by minimizing
\begin{eqnarray}\label{lossEq}
\widehat{\bbeta}_p=arg\,\, min_{\bbeta_{p}} \sumi
\rho_p\big({y_i-\bx^{\top}_i\bbeta_p}\big),
\end{eqnarray}
where $\rho_p(.)$ is as in (\ref{pdfAL}) and $\widehat{\bbeta}_p$ is
the quantile regression estimate for $\bbeta_p$ at the $p$th
quantile. The  special  case  $p=0.5$ corresponds to median
regression. As the check function is not differentiable at zero,
we cannot derive explicit solutions to the minimization problem.
Therefore, linear programming methods are commonly applied to obtain
quantile regression estimates for $\bbeta_p$.  A connection between
the minimization of the sum in (\ref{lossEq}) and the
maximum-likelihood theory is provided by the ALD; see \cite{geraci07}. It is also true that  under
the   quantile   regression model,  we   have
\begin{equation}\label{Wi}
W_i=\frac{1}{\sigma}\rho_p\big(y_i-\bx^{\top}_i\bbeta_p\big)\sim \exp(1).
\end{equation}
The  above  result is  useful  to check the model in practice, as
will be seen in the Application Section.

Now, suppose $y_1, \ldots,y_n$ are   independent observations such  as $Y_i \sim ALD (\xp^{\top}_i\bbeta_p,\sigma,p),$
$i=1,\ldots,n$. Then,  from (\ref{pdfALs})  the
log--likelihood function for $\btheta=(\bbeta_p^{\top},\sigma)^{\top} $ can be  expressed  as
\begin{equation}\label{likel}
\ell(\btheta)=\sumas \ell_i(\btheta),
\end{equation}
where
$\ell_i(\btheta)=c-\frac{3}{2}\log{\sigma}+\frac{\vartheta_p}{\tau_p^2\sigma}(y_i-\xp^{\top}_i\bbeta_p)+\log(A_i)$,
with  $c$ is a constant does not   depend on $\btheta$ and
$A_i=2\big({\frac{\delta_i}{\gamma}}\big)^{1/2}K_{1/2}(\lambda_i)=\frac{\sqrt{2\pi}}{\gamma}\exp(-\lambda_i),$
with $\delta_i=\delta(y_i)={|y_i-\xp^{\top}_i\bbeta_p|}/{\tau_p\sqrt{\sigma}})$   and $\lambda_i=\delta_i\gamma$.

Note that if we consider $\sigma$ as a nuisance parameter, then the
maximization of the likelihood in (\ref{likel}) with respect to the
parameter $\bbeta_p$ is equivalent to the minimization of the
objective function in (\ref{lossEq}). and  hence the relationship
between the check function and ALD can be used to reformulate the QR
method in the likelihood framework.

The log--likelihood function is not differentiable at zero.
Therefore,  standard  procedures  the estimation      can not   be
developed   following  the  usual way. Specifically, the  standard
errors  for    the   maximum likelihood    estimates is  not   based
on    the     genuine information matrix. To  overcome   this
problem   we   consider the  empirical information matrix as    will
be   described in the  next  Subsection.


\subsection{Parameter estimation via the EM algorithm}
In this section, we discuss  an  estimation method  for QR  based on
the EM algorithm to obtain ML estimates. Also, we consider the
method of moments (MM) estimators,which can be effectively used as
starting values in the EM algorithm.
Here,  we show how to employ the EM  algorithm
for ML estimation in QR model under the ALD. From the hierarchical
representation (\ref{hierar1})-(\ref{hierar2}), the QR model in
(\ref{QRmodel}) can be presented  as
\begin{eqnarray}
Y_i|U_i=u_i&\sim& N(\bx^{\top}_i\bbeta_p+\vartheta_p\label{repHier1}
u_i,\tau_p^2\sigma u_i),\\ U_i&\sim&
\exp(\sigma),\,\,\,\,\is,\label{repHier2}
\end{eqnarray}
where $\vartheta_p$ and $\tau_p^2$  are as in (\ref{st-ALD}). This
hierarchical representation of the QR model is  convenient to  describe
the steps  of the   EM algorithm. Let $\yp = (y_1, \ldots, y_n)$ and  $\bu = (u_1,\ldots , u_n)$ be
the observed data and   the missing data, respectively. Then,   the
complete data log-likelihood function of
$\btheta=(\bbeta^{\top}_p,\sigma)^{\top}$, given $(\yp,\bu)$,
ignoring additive constant terms, is given by
$\ell_{c}(\btheta|\by,\bu)=\sum_{i=1}^{n}\ell_{c}(\btheta|y_i,u_i)$,
where
\begin{eqnarray*}
\ell_{c}(\btheta|y_i,u_i) =-\frac{1}{2} \log( 2\pi\tau_p^2)
 -\frac{3}{2} \log(\sigma)   -\frac{1}{2}\log (u_i) -  \frac{1}{2\sigma\tau_p^2}
{u^{-1}_i}(y_i-\mathbf{x}^{\top}_i\bbeta_p-\vartheta_p u_i)^2 -
\frac{1}{\sigma} u_i,
\end{eqnarray*}
for $i=1,\ldots,n$.  In what follows the superscript $(k)$ indicates
the estimate of the related parameter  at the stage $k$ of the
algorithm. The E-step of the EM algorithm requires evaluation of the
so-called Q-function $Q(\btheta|\btheta^{(k)}) =
\textrm{E}_{\scriptsize \btheta^{(k)}}[\ell_{c}(\btheta|\y,\bu)|\by, \btheta^{(k)}]$, where
$\textrm{E}_{\scriptsize \btheta^{(k)}}[.]$ means that the
expectation is being effected using $\btheta^{(k)}$ for $\btheta$.
Observe that the expression of the Q-function is completely
determined by the knowledge of the expectations
\begin{eqnarray} \label{weith}
{\cal E}_{ s i}(\btheta^{(k)}) =  \textrm{E}[U^s_i |y_i, \btheta^{(k)}],\,\,\, s=-1,1,
\end{eqnarray}
that  are obtained of   properties of   the $GIG(0.5, a, b)$ distribution. Let   $\bxi^{(k)} _{s}=\big( {\cal
E}_{s1}(\btheta^{(k)}), \ldots, {\cal E}_{sn}(\btheta^{(k)})
\big)^{\top}$ be the vector  that  contains  all   quantities   defined
in (\ref{weith}). Thus, dropping unimportant constants,    the
Q-function can be written in a synthetic form as
$Q(\btheta|\widehat{\btheta})=\sum_{i=1}^{n}Q_i(\btheta|\widehat{\btheta})$,
where {\small{
\begin{eqnarray} \label{eqn qfunction}
Q_i(\btheta|\widehat{\btheta}) =
-\frac{3}{2}\log\sigma-\frac{1}{2\sigma\tau_p^2}
\left[ {\cal
E}_{-1i}(\btheta^{(k)})(y_i-\mathbf{x}^{\top}_i\bbeta_p)^2-2(y_i-\mathbf{x}^{\top}_i\bbeta_p)\vartheta_p+\frac{1}{4}{\cal E}_{1
i}(\btheta^{(k)})\tau_p^4\right]. \, \, \, \, \,
\end{eqnarray}}}
This  quite useful expression to implement the M-step, which consists of  maximizing it over $\btheta$. So the  EM algorithm   can be summarized as  follows:\\
\noindent \emph{E-step}: Given $\btheta=\btheta^{(k)}$,  compute ${\cal E}_{si}(\btheta^{(k)})$ through  of the   relation
\begin{equation}\label{deltai}
{\cal E}_{ si}(\btheta^{(k)}) =E[U^s_i|y_i,\btheta^{(k)}]=\left(\frac{\delta^{(k)}_i}{\gamma^{(k)}}\right)^{s}\frac{K_{1/2+s}\big(\lambda^{(k)}_i\big)}{K_{1/2}\big(\lambda^{(k)}_i\big)}, s=-1,1,
\end{equation}
where
$\delta^{(k)}_i=\frac{|y_i-\xp^{\top}_i\bbeta^{(k)}_p|}{\tau_p\sqrt{\sigma^{(k)}}}$,  $\gamma^{(k)}=\frac{\tau_p}{2\sqrt{\sigma^{(k)}}}$ and
$\lambda^{(k)}_i={\delta^{(k)}_i\gamma^{(k)}}$;\\
\noindent{\em  M-step}: Update ${\btheta}^{(k)}$ by
maximizing $Q(\btheta|\btheta^{(k)})$  over $\btheta$,
which leads to the following expressions
{\small{
\begin{eqnarray*}
{\bbeta}^{(k+1)}_p&=&\left(\X^{\top} D(\bxi^{(k)} _{-1})\X \right)^{-1}\X^{\top}\big(D(\bxi^{(k)} _{-1})\Y-\vartheta_p {\bf 1}_n\big),\, \, \\
{\sigma}^{(k+1)}&=&\frac{1}{3n\tau^2_p}\Big[Q(\bbeta^{(k+1)}, \bxi_{-1}^{(k)})-2{\bf 1}^{\top}_n(\Y-\X
\bbeta^{(k+1)})\vartheta_p+\frac{\tau_p^4}{4}{\bf
1}^{\top}_n\bxi^{(k)} _{1}\Big],
\end{eqnarray*}}}
where $D(\mathbf{a})$ denotes the diagonal matrix, with the diagonal
elements given by $\mathbf{a}=(a_1,\ldots,a_p)^{\top}$ and
$Q(\bbeta, \bxi_{-1})= (\Y-\X \bbeta)^{\top}D(\bxi_{-1}) (\Y-\X\bbeta)$. A similar  expression  for  $\bbeta^{(k+1)}_p$ is  obtained in  \cite{tian13}.
This process is iterated until some distance involving two
successive evaluations of the actual log-likelihood $\ell(\btheta)$,
like $||\ell({\btheta}^{(k+1)})-\ell({\btheta}^{(k)})||$ or
$||\ell({\btheta}^{(k+1)})/\ell({\btheta}^{(k)})-1||$, is small
enough. This algorithm is implemented as part of the R package
\verb"ALDqr()",  which can be downloaded at not cost from the
repository CRAN. Furthermore, following  the  results  given in
\cite{Yu2005}, the MM estimators for $\bbeta_p$ and $\sigma$ are
solutions of the following equations:
\begin{eqnarray}\label{ab;m}
\widehat{\bbeta}_{pM}=\big(\X^{\top}
\X\big)^{-1}\X^{\top}\big(\Y-\widehat{\sigma}_M \vartheta_p{\bf
1}_n\big)\, \, \, \,{\rm and} \, \,
\widehat{\sigma}_{M}=\displaystyle  \frac{1}{n}\sum_{i=1}^n
\rho_p\big(y_i-\bx_i^{\top}\widehat{\bbeta}_{pM}\big),
\end{eqnarray}
where $\vartheta_p$  is as (\ref{st-ALD}). Note that the
MM estimators do not have explicit closed form and numerical
procedures are needed to solve these non-linear equations. They can be used  as  initial values in the iterative procedure for computing the ML estimates  based on the EM-algorithm.
Standard   errors  for    the   maximum  likelihood    estimates is  based  on
the  empirical information matrix, that  according to \cite{meilijson89} formula, is defined as
\begin{eqnarray}\label{Imatrix }
\mathbf{L}(\btheta)=\sum_{j=1}^{n}\textbf{s}(y_j|\btheta)\textbf{s}^\top(y_j|\btheta)-n^{-1}\textbf{S}(y_j|\btheta)\textbf{S}^\top(y_j|\btheta),
\end{eqnarray}
where $\textbf{S}(y_j|\btheta)=\sum_{j=1}^{n}\textbf{s}(y_j|\btheta)$. It is noted from the result of \cite{louis82} that the individual score can be determined as
$\textbf{s}(y_j|\btheta) ={\partial \log f(y_j|\btheta)}/{\partial \btheta} = E\Big({\partial \ell_{c_j}(\btheta|y_j, u_i)}/{\partial \btheta} |  y_j,\btheta\Big)$.  Asymptotic
confidence intervals and tests of the parameters at the $p$th
level can be obtained assuming that the ML estimator
$\widehat{\btheta}$ has approximately a normal  multivariate distribution. \\

From  the EM algorithm, we  can see that $ {\cal E}_{-1 i}(\btheta^{(k)})$ is  inversely proportional to $d_i=|y_i-\xp^{\top}_i\bbeta^{(k)}_p|/\sigma$. Hence, $u_i(\btheta^{(k)})= {\cal E}_{-1 i}(\btheta^{(k)})$   can be
interpreted as a type of weight for the $i$th case in the estimates
of $\bbeta^{(k)}_p$, which   tends to be small for   outlying
observations. The behavior of    these weights can be used as tools for
identifying outlying observations   as  well as  for  showing that    we   are   considering a robust approach,  as will be seen in Sections 4 and 5.

\section{Case-deletion measures} \label{Sec Diagnostic}
Case-deletion is a classical  approach to study the effects of
dropping the $i$th case from the data set. Let $\y_{c}=(\y,\bu)$ be  the augmented data set,   and  a quantity
with a subscript ``$[i]$''  denotes the original one with the $i$th
observation deleted. Thus, The  complete-data log-likelihood function
based on  the data with the $i$th case deleted will be denoted by
$\ell_{c}(\btheta|\y_{c[i]})$.  Let
$\widehat{\btheta}_{[i]}=(\widehat{\bbeta}^{\top}_{p[i]},
\widehat{\sigma^2}_{[i]})^{\top}$ be the  maximizer of the function
{{ $Q_{[i]}(\btheta|\widehat{\btheta})=
\textrm{E}_{\scriptsize{\widehat{\btheta}}}\left[\ell_{c}(\btheta|\Y_{c[i]})|\y
\right] $}},  where $\widehat{\btheta}=(\widehat{\bbeta}^{\top},
\widehat{\sigma^2})^{\top}$ is the ML estimate of $\btheta$. To
assess the influence of the $i$th case on $\widehat{\btheta}$, we
compare the difference between $\widehat{\btheta}_{[i]}$  and
$\widehat{\btheta}$. If the deletion of a case seriously influences
the estimates, more attention needs  to  be paid to that case.
Hence, if $\widehat{\btheta}_{[i]}$ is far from  $\widehat{\btheta}$
in some sense, then the $i$th case is regarded as influential. As
$\widehat{\btheta}_{[i]}$ is needed for every case, the required
computational  effort can be quite heavy, especially when the sample
size  is large. Hence, To calculate the case-deletion estimate $\widehat{\btheta}^1_{[i]}$ of $\btheta$,  \citep[see][]{zhu2001case} proposed the following one-step  approximation based on the Q-function,
\begin{eqnarray}\label{theta1}
\widehat{\btheta}^1_{[i]}= \widehat{\btheta}+ \big\{
-\ddot{Q}(\widehat{\btheta}|\widehat{\btheta})\big\}^{-1}
\dot{Q}_{[i]}(\widehat{\btheta}|\widehat{\btheta}),
\end{eqnarray}
where
\begin{eqnarray} \label{eqn Hessian Matrix and Grad}
\ddot{Q}(\widehat{\btheta}|\widehat{\btheta})=\displaystyle\frac{\partial^2
Q(\btheta|\widehat{\btheta})}{\partial\mbox{\btheta}\partial{\btheta}^{\top}}
\big\vert_{\btheta=\widehat{\btheta}}   \,\,\,\, \textrm{and}
\,\,\,\,
  \dot{Q}_{[i]}(\widehat{\btheta}|\widehat{\btheta})=\displaystyle\frac{\partial{{Q}_{[i]}(\btheta|\widehat{\btheta})}}{\partial{\btheta}}\big\vert_{\btheta=\widehat{\btheta}},
\end{eqnarray}
are  the Hessian matrix  and the gradient vector evaluated at
$\widehat{\btheta}$, respectively. The Hessian matrix
is an essential element in the method  developed by \cite{zhu2001case} to obtain the measures for
case-deletion diagnosis. For    developing  the case-deletion measures, we  have   to   obtain   the  elements in   (\ref{theta1}), $\dot{Q}_{[i]}(\widehat{\btheta}|\widehat{\btheta})$ and  $\ddot{Q}(\widehat{\btheta}|\widehat{\btheta})$.  These formulas  can be obtained quite easily from (\ref{eqn qfunction}):

\begin{enumerate}
\item[1.]  The  components  of $\dot{Q}_{[i]}(\widehat{\btheta}|\widehat{\btheta})$ are
{\small{
\begin{eqnarray*}
\dot{Q}_{[i]\bbeta}(\widehat{\btheta}|\widehat{\btheta})&=& \frac{  \partial{Q_{[i]}({\btheta}|\widehat{\btheta})}}{\partial{\bbeta}}\big\vert_{\btheta=\widehat{\btheta}} ={{\frac{1}{\widehat{\sigma}}}} E_{1[i]}
\end{eqnarray*}
and
\begin{eqnarray*}
\,\, \,
\dot{Q}_{[i]\sigma}(\widehat{\btheta}|\widehat{\btheta})=
\frac{\partial{Q_{[i]}({\btheta}|\widehat{\btheta})}}{
\partial{\sigma}}\big\vert_{\btheta=\widehat{\btheta}}=-\frac{1}{2\widehat{\sigma^2}}
E_{2[i]},    \label{sigma}
\end{eqnarray*}
where
\begin{eqnarray}
 E_{1[i]} & = & \frac{1}{\tau_{p}^{2}} \sum_{j\neq i} \left[{\cal E}_{-1j}(\widehat{\btheta}^{(k)})(y_{j}-\mathbf{x}_{j}^{\top}\widehat{\bbeta})\mathbf{x}_{j}-\mathbf{x}_{j}\vartheta_{p} \right] \,\,\,\, \textrm{ and} \label{eqn E1i}  \\
E_{2[i]} & = & \sum_{j\neq i} \left[3
\widehat{\sigma}-\frac{1}{\tau_{p}^{2}}  {\cal
E}_{-1j}(\widehat{\btheta}^{(k)})(y_{j}-\mathbf{x}_{j}^{\top}\widehat{\bbeta_p})^{2}-2(y_{j}-\mathbf{x}_{j}^{\top}\widehat{\bbeta_p})\vartheta_{p}
+ \frac{1}{4}{\cal
E}_{1j}(\widehat{\btheta}^{(k)})\tau_{p}^{4}\right].
\label{eqn E2i}
\end{eqnarray}}}
\item[2.]  The elements   of  the  second order partial derivatives of $Q(\btheta|\widehat{\btheta})$ evaluated at $\widehat{\btheta}$ are
{{
\begin{eqnarray*}
\ddot{Q}_{\beta}(\widehat{\btheta}|\widehat{\btheta})& = &  -\frac{1}{\widehat{\sigma}\tau_p^{2}} \X^{\top}D\big(\bxi _{-1}^{(k)}\big)\X,\nonumber\\
\ddot{Q}_{\sigma}(\widehat{\btheta}|\widehat{\btheta})\}&
=& \frac{3}{4\widehat{\sigma^2}} - \frac{1}{2\widehat{\sigma^3}\tau^2_p}\Big[Q\big(\bbeta,\bxi^{(k)} _{-1}\big)-2{\bf 1}^{\top}_n(\Y-\X
\bbeta)\vartheta_p+\frac{\tau_p^4}{4}{\bf
1}^{\top}_n\bxi^{(k)} _{1}\Big]
\end{eqnarray*}}}
and  $\ddot{Q}_{\beta \sigma}(\widehat{\btheta}|\widehat{\btheta})\}
= {\bf 0}$.
\end{enumerate}
In   the   following result,  we will  obtain the
one-step approximation of
$\widehat{\btheta}_{[i]}=(\widehat{\bbeta}^{\top}_{p[i]},
\widehat{\sigma}_{[i]})^{\top}$, $i=1,\ldots,n$ based on
(\ref{theta1}), viz., the relationships between  the parameter
estimates for the full data set and the data with the $i$th case
deleted.
\begin{theorem} \label{the;1}
For the QR model defined in (\ref{repHier1}) and (\ref{repHier2}), the
relationships between the parameter estimates for full data set and
the data with the $i$th case deleted are  as follows:
\begin{eqnarray*}\label{aprox1}
\widehat{\bbeta}^1_{p[i]}&=& \widehat{\bbeta}_p+  \tau_p^{2} \big(\X^{\top}D\big(\widehat{\bxi}_{-1}\big)\X\big)^{-1} \textbf{E}_{1[i]}\,\,\, \, \,{\rm and }\,\, \, \, \,
\widehat{\sigma^2}^1_{[i]}= \widehat{\sigma^2} - \frac{1}{2\widehat{\sigma^2}}\Big(\ddot{Q}_{\sigma}(\widehat{\btheta}|\widehat{\btheta})\Big)^{-1} E_{2[i]},
\end{eqnarray*}
where  $\textbf{E}_{1[i]}$ and   $E_{2[i]}$ are as in (\ref{eqn E1i}) and  (\ref{eqn E2i}), respectively.
\end{theorem}
To asses the influence of the $i$th case on the ML estimate
$\widehat{\btheta}$, we  compare $\widehat{\btheta}_{[i]}$
and $\widehat{\btheta}$ based on  metrics,  proposed by \cite{zhu2001case},  for measuring the
distance between $\widehat{\btheta}_{[i]}$ and  $\widehat{\btheta}$.  For  that, we consider here the following;

\begin{enumerate}
\item  {\it Generalized Cook distance}:
\begin{equation}\label{GCD}
GD_i=(\widehat{\btheta}_{[i]}-\widehat{\btheta})^{\top}\big\{ -\ddot{
Q}(\widehat{\btheta}|\widehat{\btheta})\big\}(\widehat{\btheta}_{[i]}-\widehat{\btheta}),
\quad \is.
\end{equation}
Upon  substituting  (\ref{theta1}) into  (\ref{GCD}), we obtain the
approximation
$$GD^1_i=\dot{Q}_{[i]}(\widehat{\btheta}|\widehat{\btheta})^{\top}\big\{-\ddot{Q}(\widehat{\btheta}|\widehat{\btheta})\big\}^{-1}
\dot{Q}_{[i]}(\widehat{\btheta}|\widehat{\btheta}), \quad \is.$$
As $\ddot{Q}(\widehat{\btheta}|\widehat{\btheta})$ is  a  diagonal matrix, one  can    obtain   easily a  type of  Generalized Cook distance for  parameters $\bbeta$ and $\sigma$, respectively,  as  follows
$$GD^1_i(\bbeta)=\dot{Q}_{[i]\bbeta}(\widehat{\btheta}|\widehat{\btheta})^{\top}\big\{-\ddot{Q}_{\beta}(\widehat{\btheta}|\widehat{\btheta})\big\}^{-1}
\dot{Q}_{[i]\bbeta}(\widehat{\btheta}|\widehat{\btheta}), \quad \is.$$
$$GD^1_i(\sigma)=\dot{Q}_{[i]\sigma}(\widehat{\btheta}|\widehat{\btheta})^{\top}\big\{-\ddot{Q}_{\sigma}(\widehat{\btheta}|\widehat{\btheta})\big\}^{-1}
\dot{Q}_{[i]\sigma}(\widehat{\btheta}|\widehat{\btheta}), \quad \is.$$

\item  {\it  Q-distance}:  This  measure of the influence of the $i$th case is based  on  the
$Q$-distance function, similar to the likelihood distance $LD_i$
\citep{cook82},  defined as
\begin{equation}\label{QD}
QD_i=2\big\{Q(\widehat{\btheta}|\widehat{\btheta})-Q(\widehat{\btheta}_{[i]}|\widehat{\btheta})\big\}.
\end{equation}
We can calculate an approximation of the  likelihood displacement
$QD_i$ by  substituting  (\ref{theta1}) into  (\ref{QD}), resulting
in  the following approximation $QD^{1}_i$ of $QD_i$:
\begin{equation*}\label{QD1}
QD^1_i=2\big\{Q(\widehat{\btheta}|\widehat{\btheta})-Q(\widehat{\btheta}^1_{[i]}|\widehat{\btheta})\big\}.
\end{equation*}
\end{enumerate}

\section{Application}\label{sec application}

\begin{figure}[!t]
\begin{center}
\includegraphics[scale=0.35]{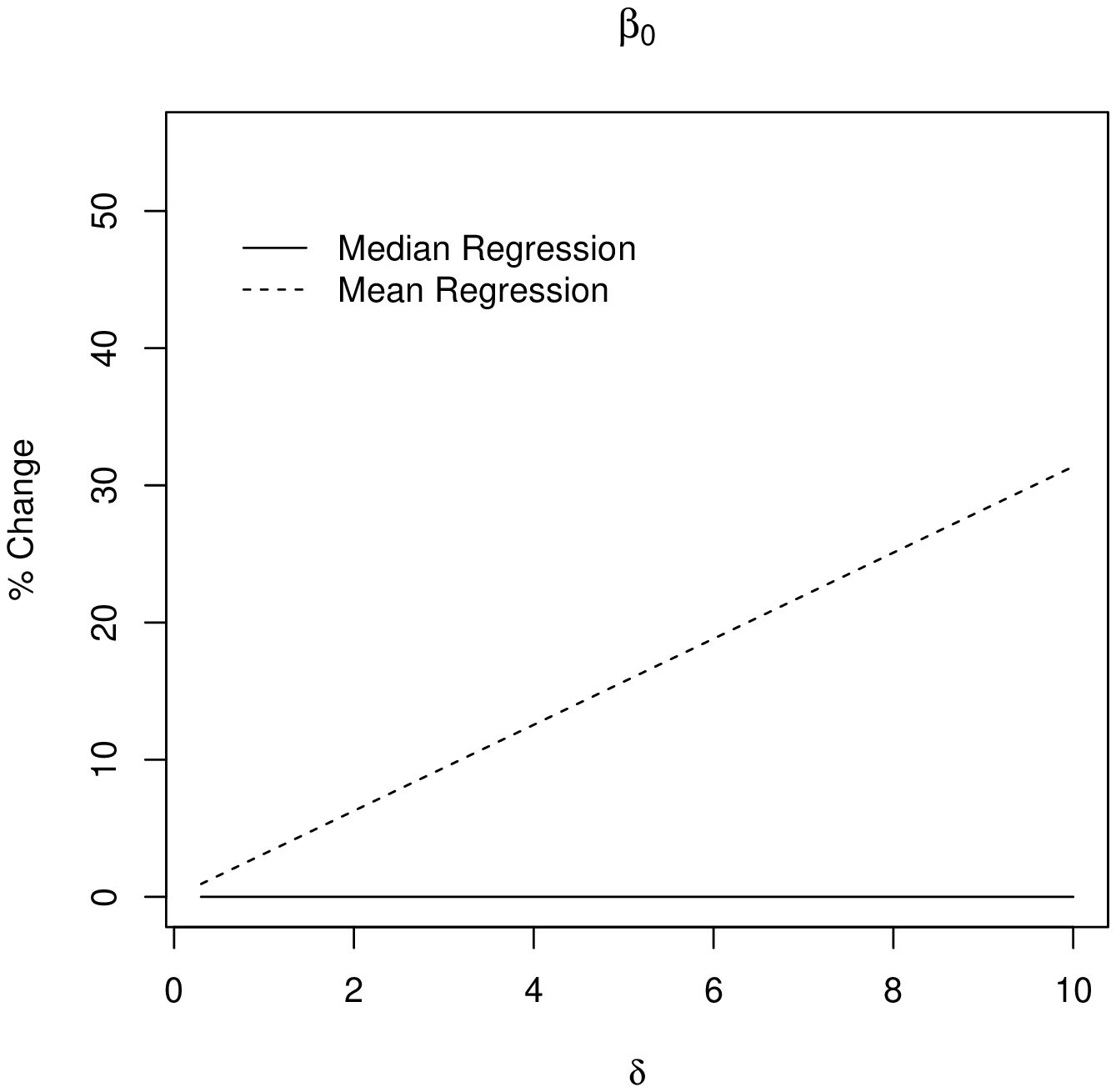}~\includegraphics[scale=0.35]{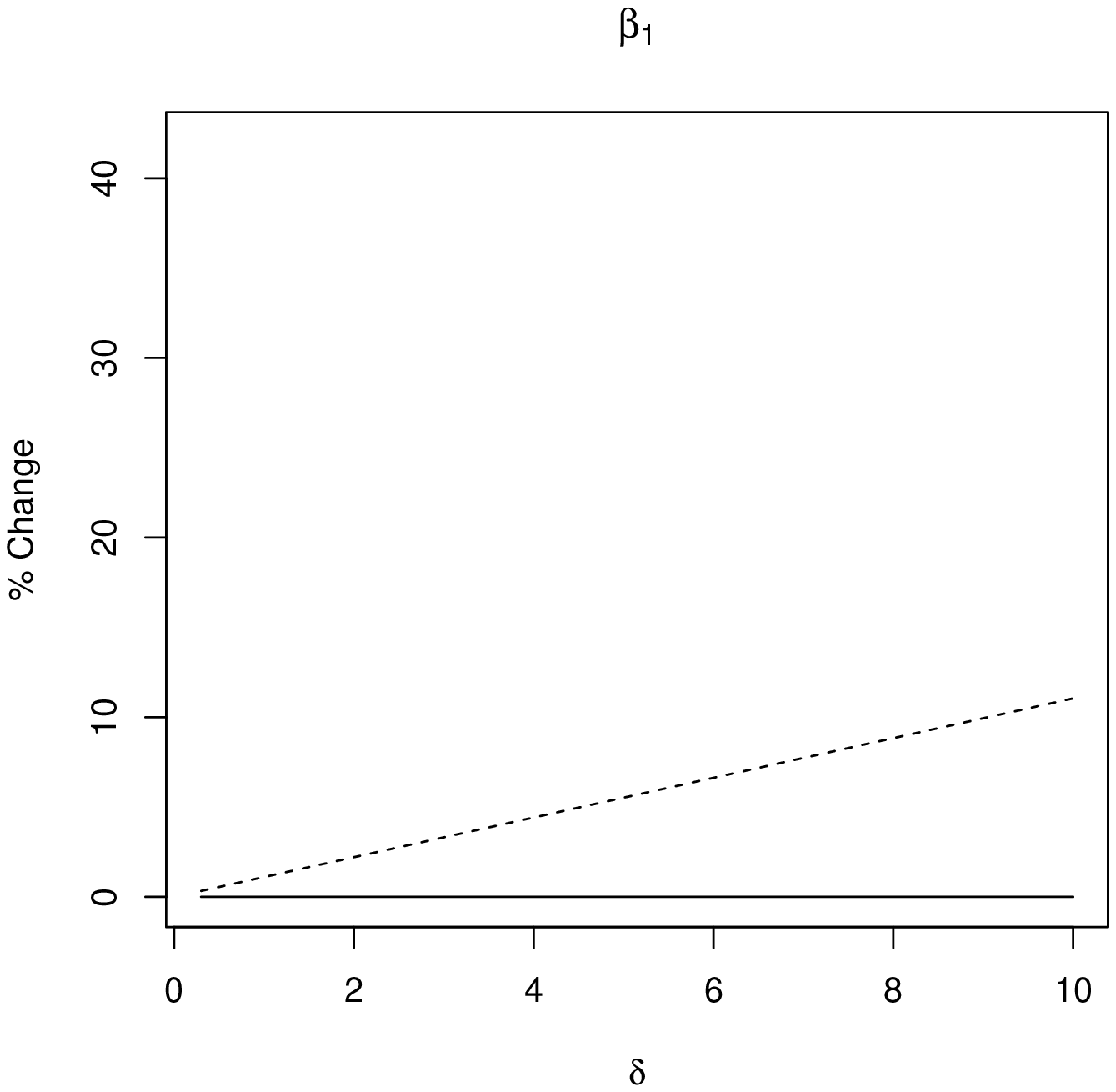}~\includegraphics[scale=0.35]{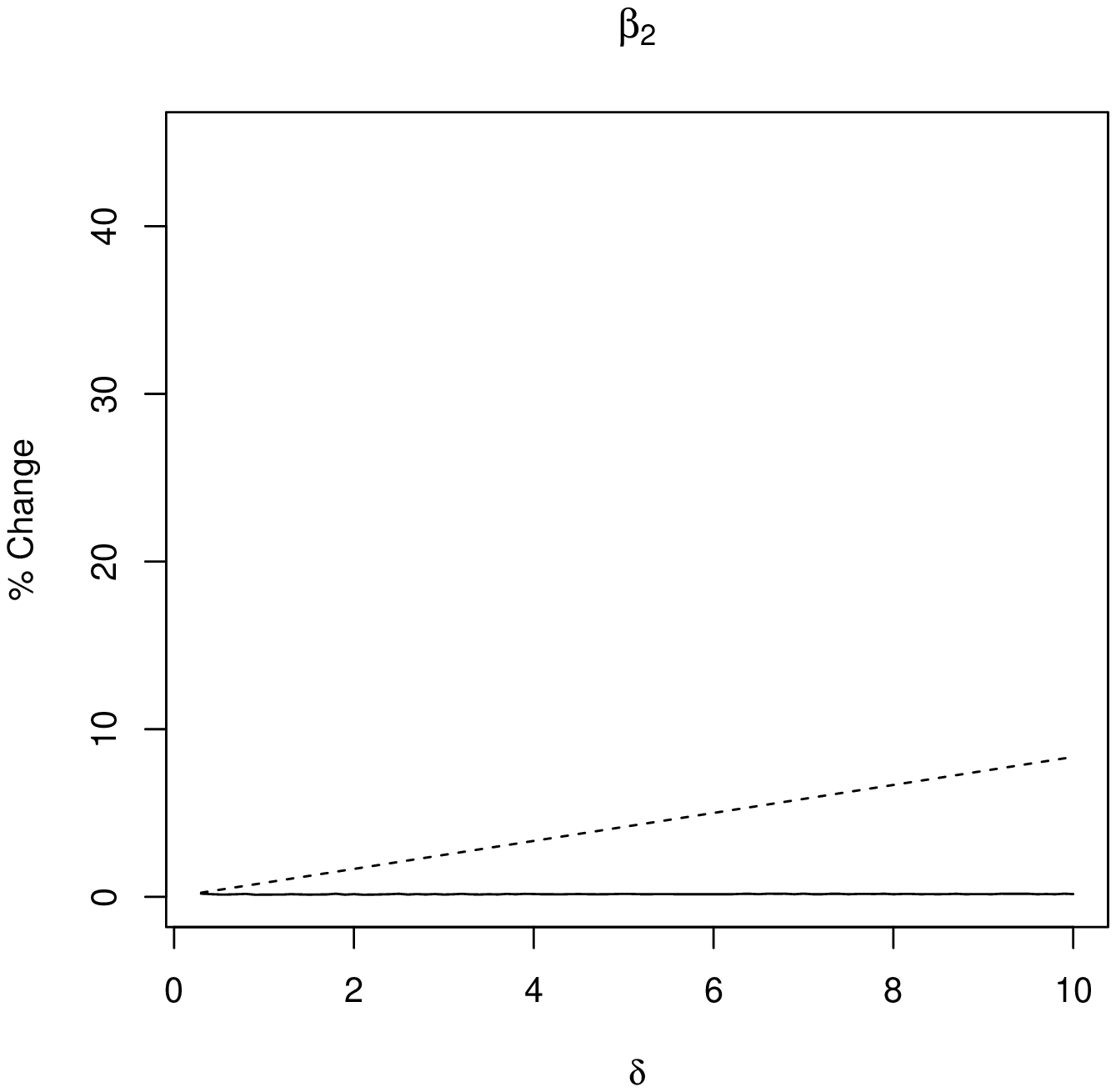}
\caption{Percentage of change in the estimation of $\beta_0$,
$\beta_1$ and $\beta_2$ in comparison with the true value, for
median $(p=0.5)$ and mean regression, for different contaminations
$\delta$. \label{fig:change}}
\end{center}
\end{figure}

\begin{figure}[!tb]
\begin{center}
\includegraphics[scale=0.55]{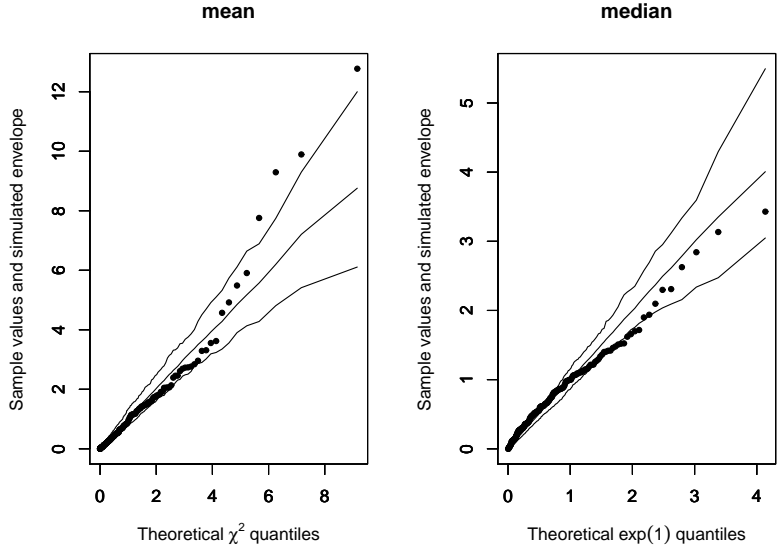}~
\caption{AIS data: Q--Q plots and simulated envelopes for mean and
median regression.\label{fig:2c}}
\end{center}
\end{figure}

\begin{figure}[!t]
\begin{center}
\includegraphics[scale=0.47]{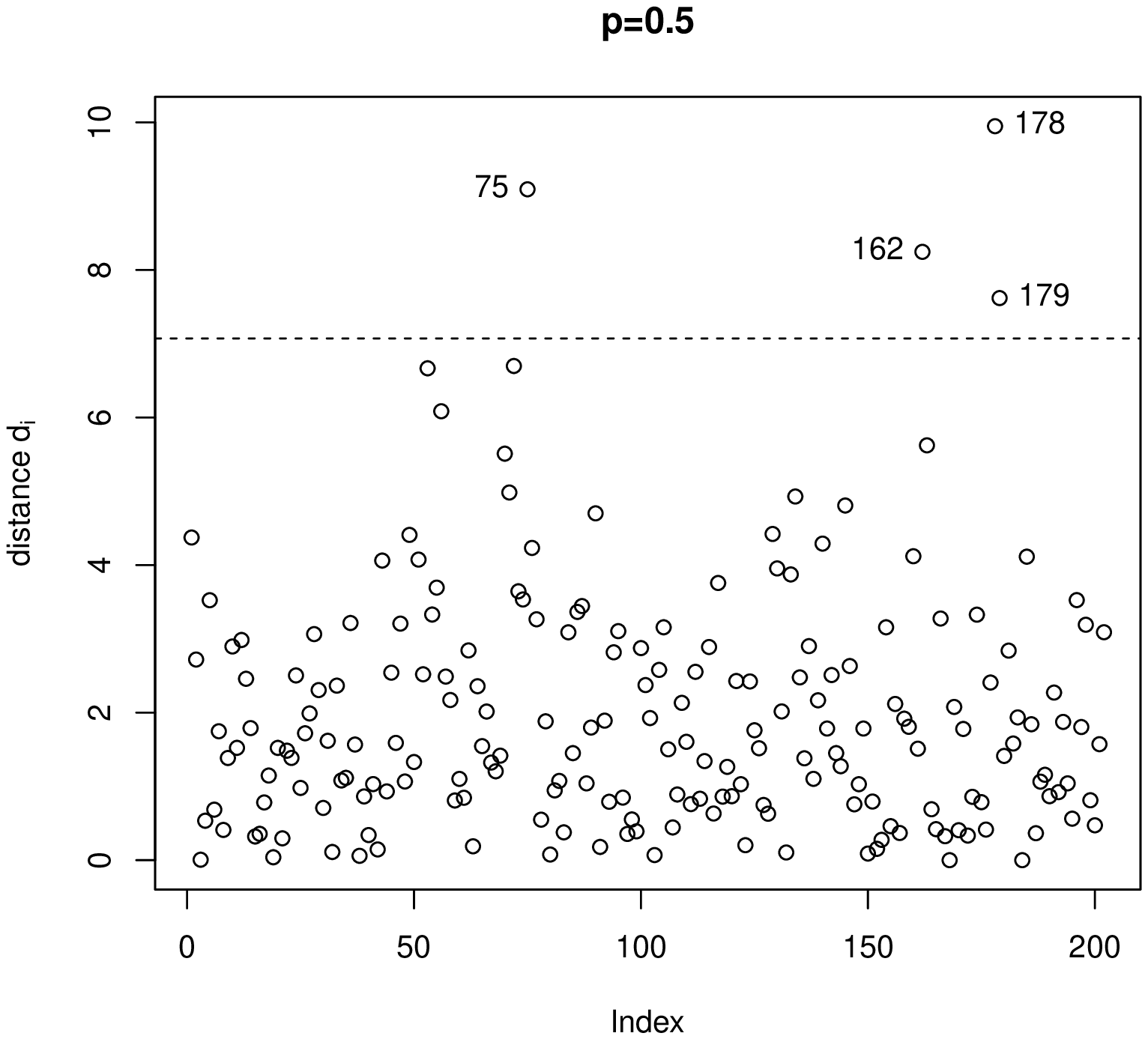}~\includegraphics[scale=0.5]{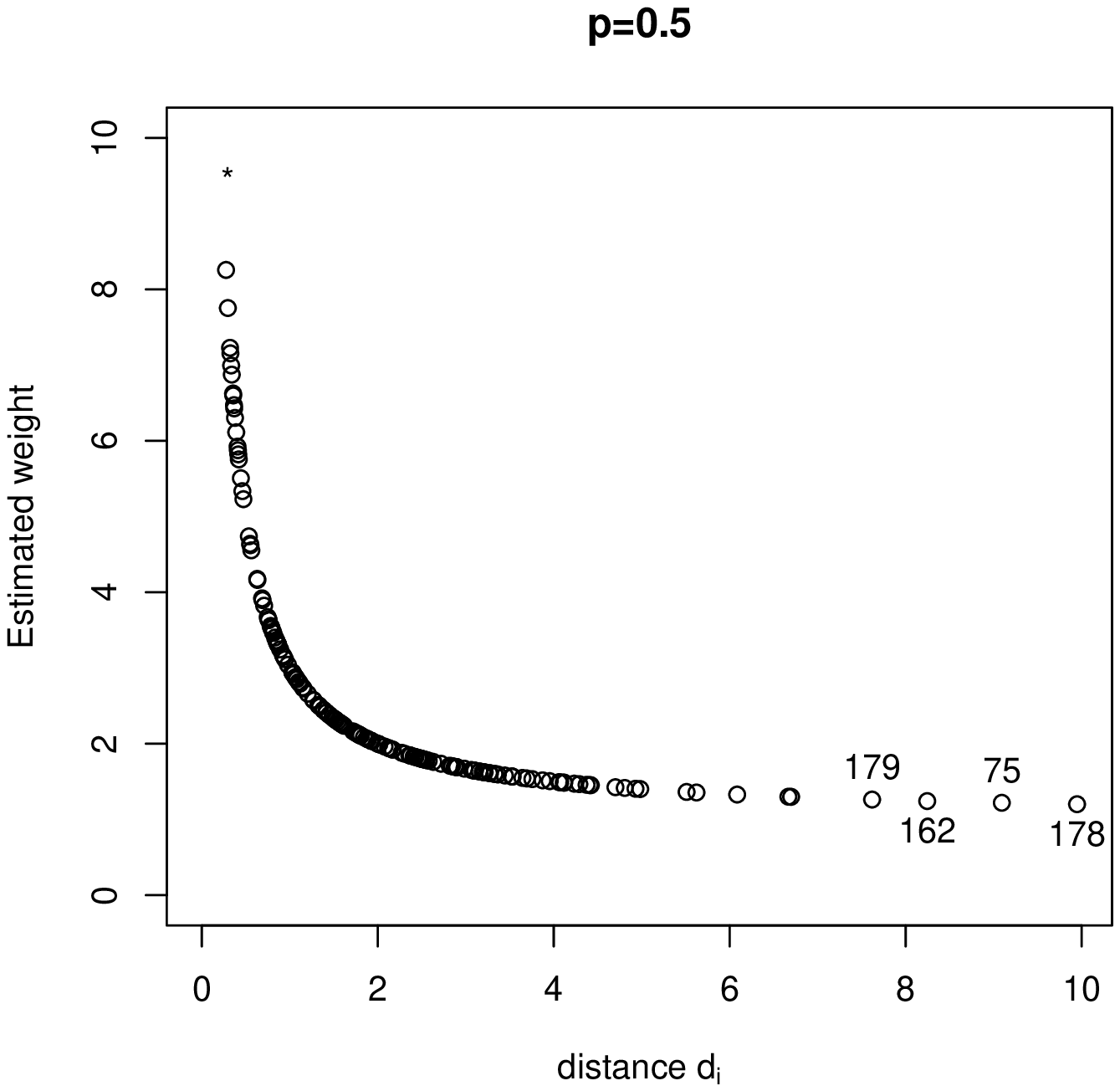}
\caption{AIS data: Index plot of the distance  $d_i$ and the
estimated weights $u_i$.\label{fig:mahal}}
\end{center}
\end{figure}

We illustrate the proposed methods by applying them to the
Australian Institute of Sport (AIS) data, analyzed by Cook and
Weisberg (1994) in a normal regression setting. The data set
consists of several variables measured in $n=202$ athletes (102
males and 100 females). Here, we focus on body mass index (BMI),
which is assumed to be explained by  lean body mass (LBM) and gender
(SEX). Thus, we consider the following QR model:
$$BMI_i=\beta_0+\beta_1 LBM_i+\beta_2 SEX_i+\epsilon_i,\,\,\,\,\,i=1,\ldots,202,$$
where $\epsilon_i$ is a zero $p$ quantile. This model can be fitted
in the R software by using the package \verb"quantreg()", where one
can arbitrarily use the BR or the LPQR algorithms. In order to
compare with our proposed EM algorithm, we carry out quantile
regression at three different quantiles, namely $p= \{0.1, 0.5,
0.9\}$ by using the ALD distribution as described in Section 2. The
ML estimates and associated standard errors were obtained by using
the EM algorithm and the observed information matrix described in
Subsections 2.3, respectively. Table \ref{table.application}
compares the results of our EM, BR and the LPQR estimates under the
three selected quantiles. The standard error of the LPQR estimates
are not provided in the R package \verb"quantreg()" and are not
shown in Table \ref{table.application}. From this table we can see
that estimates under the three methods only exhibit slight
differences, as expected. However, the standard errors of our EM
estimates are smaller than those via the BR algorithm. This suggests
that the EM algorithm seems to produce more accurate estimates of
the regression parameters at the $p$th level.
\begin{table}[ht!]
\centering
\caption{AIS data. Results of the parameter estimation via EM, Barrodale and Roberts (BR) and  Lasso Penalized Quantile Regression (LPQR) algorithms for three selected quantiles.}
{\small{
\begin{tabular}{ccccccc}

       \hline\hline
    \multicolumn{2}{c}{} & \multicolumn{2}{c}{EM} &  \multicolumn{2}{c}{BR}   &  \multicolumn{1}{c}{LPQR} \\

        $p$ & Parameter   & MLE        &     SE   & Estimative &  SE    &Estimative \\
        \hline
        0.1 & $\beta_{0}$ & 9.3913     &   0.7196 & 9.3915     & 1.2631 & 9.8573     \\
            & $\beta_{1}$ & 0.1705     &   0.0091 & 0.1705     & 0.0160 & 0.1647     \\
            & $\beta_{2}$ & 0.8312     &   0.2729 & 0.8209     & 0.4432 & 0.6684     \\
            & $\sigma$    & 0.2617     &   0.0252 & 1.0991     & ------ & 1.0959    \\
        \hline
        0.5 & $\beta_{0}$ & 7.6480     &   0.8717 & 7.6480     & 1.1120 & 7.6480    \\
            & $\beta_{1}$ & 0.2160     &   0.0116 & 0.2160     & 0.0159 & 0.2160    \\
            & $\beta_{2}$ & 2.2499     &   0.3009 & 2.2226     & 0.4032 & 2.2226    \\
            & $\sigma$    & 0.6894     &   0.0590 & 0.6894     & ------ & 0.6894    \\
        \hline
        0.9 & $\beta_{0}$ & 5.8000     &   0.5887 & 5.8000     & 1.6461 & 6.0292     \\
            & $\beta_{1}$ & 0.2700     &   0.0084 & 0.2700     & 0.0256 & 0.2678     \\
            & $\beta_{2}$ & 3.9596     &   0.1937 & 3.9658     & 0.6203 & 3.8271     \\
            & $\sigma$    & 0.3391     &   0.0258 & 1.2677     & ------ & 1.2767     \\
      \hline\hline

\end{tabular}
}}
  \label{table.application}%
\end{table}%
\begin{figure}[!t]
\begin{center}
\includegraphics[scale=0.7]{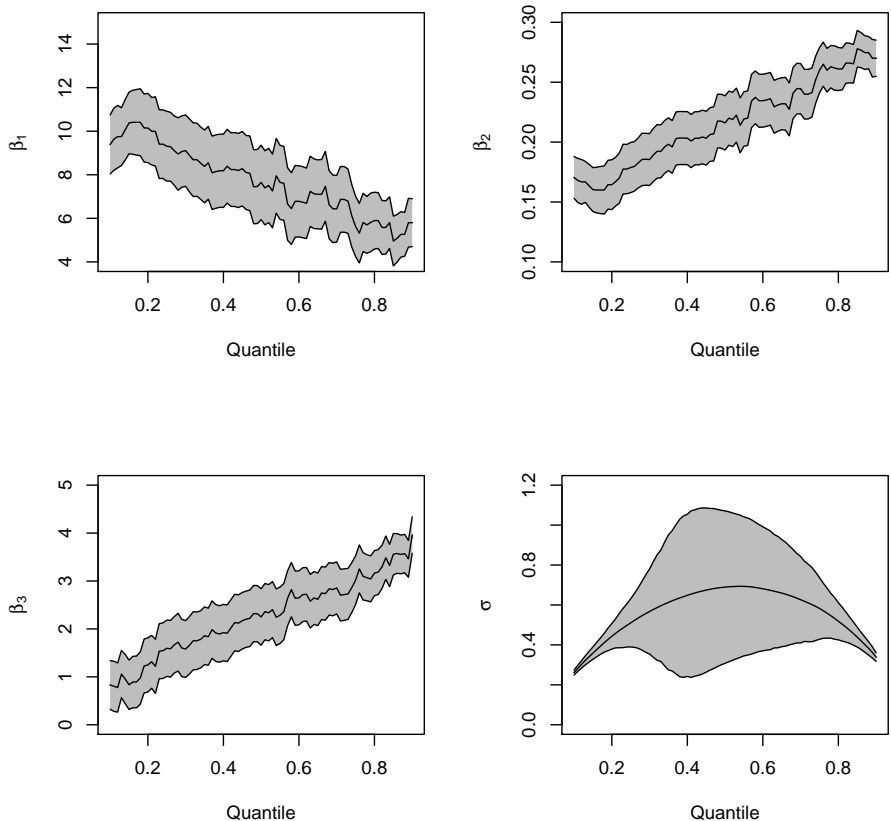}~
\caption{AIS data: ML estimates  and $95\%$ confidence intervals for
various values of $p$. \label{fig:2b}}
\end{center}
\end{figure}
To obtain a more complete picture of the effects, a series of QR
models over the grid $p=\{0.1, 0.15,\ldots, 0.95\}$ is estimated.
Figure  \ref{fig:2b} gives a graphical summary of this analysis. The
shaded area depicts the $95\%$ confidence interval from all the
parameters. From Figure \ref{fig:2b} we can observe some interesting
evidences which cannot be detected by mean regression. For example,
the effect of the two variables (LBM and gender) become stronger for
the higher conditional quantiles,  indicating that the BMI are
positively correlated with the quantiles. The robustness of the
median regression $(p=0.5)$ can be assessed by considering the
influence of a single outlying observation on the EM estimate of
$\btheta$. In particular, we can assess how much the EM estimate of
$\btheta$ is influenced by a change of $\delta$ units in a single
observation $y_{i}$. Replacing $y_{i}$ by
$y_{i}(\delta)=y_{i}+\delta sd(\mathbf{y})$, where $sd(.)$ denotes
the standard deviation. Let $\widehat{\beta}_{j}(\delta)$ be the EM
estimates of $\beta_j$ after contamination, $j=1,2,3$. We are
particularly interested in the relative changes
$|(\widehat{\beta}_{j}(\delta)-\widehat{\beta}_{j})/\widehat{\beta}_{j}|$.
In this study we contaminated the observation corresponding to
individual $\left\{\#146\right\}$ and for $\delta$ between 0 and 10.
Figure \ref{fig:change} displays  the results of the relative
changes of the estimates for different values of $\delta$. As
expected, the estimates from the median regression model are less
affected by variations on $\delta$ than those of the mean
regression. Moreover,  Figure \ref{fig:2c} shows  the Q-Q plot and envelopes  for mean and
median regression,  which  are  obtained based on the distribution of
$W_i$,  given in (\ref{Wi}), that  follows  $\exp(1)$ distribution. The lines in these figures represent the 5th
percentile, the mean and the $95$th percentile of $100$ simulated
points for each observation. These figures clearly show that the
median regression
distribution provides a better-fit  than the standard mean regression to the AIS data set.\\
As discussed at  the  end  of  Section 2.3 the estimated distance
$\widehat{d}_i=|y_i-\xp^{\top}_i\widehat{\bbeta}_p|/\widehat{\sigma}$
can be used efficiently as a  measure  to identify possible outlying
observations. Figure \ref{fig:mahal}(left panel) displays the index
plot of the distance $d_i$ for the median regression model
$(p=0.5)$. We see from this figure that observations {\#75, \#162,
\#178 and \#179} appear as possible outliers. From the EM-algorithm,
the estimated weights $u_i(\widehat{\btheta})={\cal E}_{
si}(\widehat{\btheta})$ for these observations are the smallest ones
(see right panel in Figure \ref{fig:mahal}), confirming the
robustness aspects of the maximum likelihood estimates against
outlying observations of the QR models. Thus, larger $d_i$ implies a
smaller $u_i(\widehat{\btheta})$, and the estimation of $\btheta$
tends to give smaller weight to outlying observations in the sense
of the distance $d_i$.

\begin{figure}[!t]
\begin{center}
\includegraphics[scale=0.5]{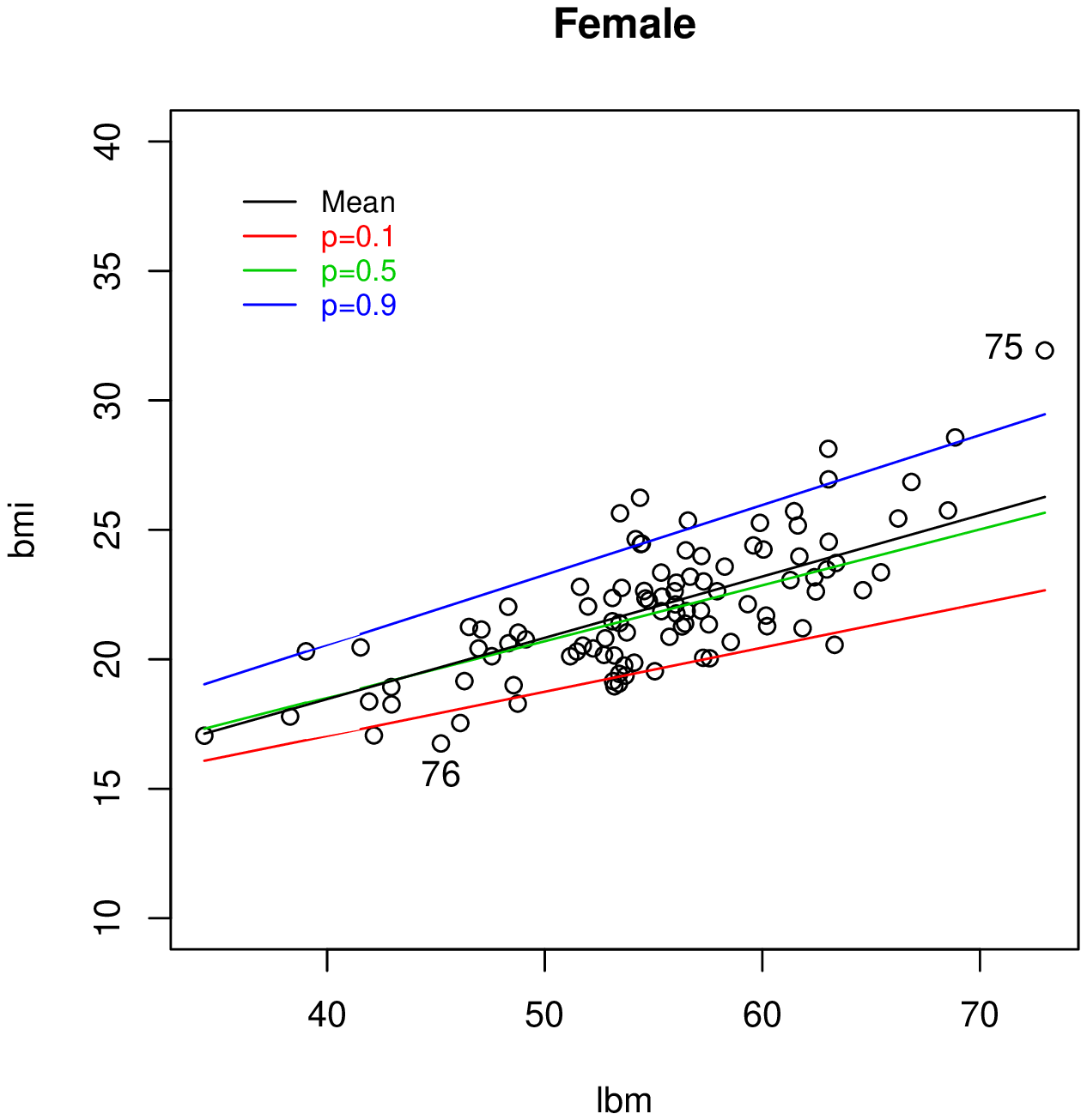}~\includegraphics[scale=.5]{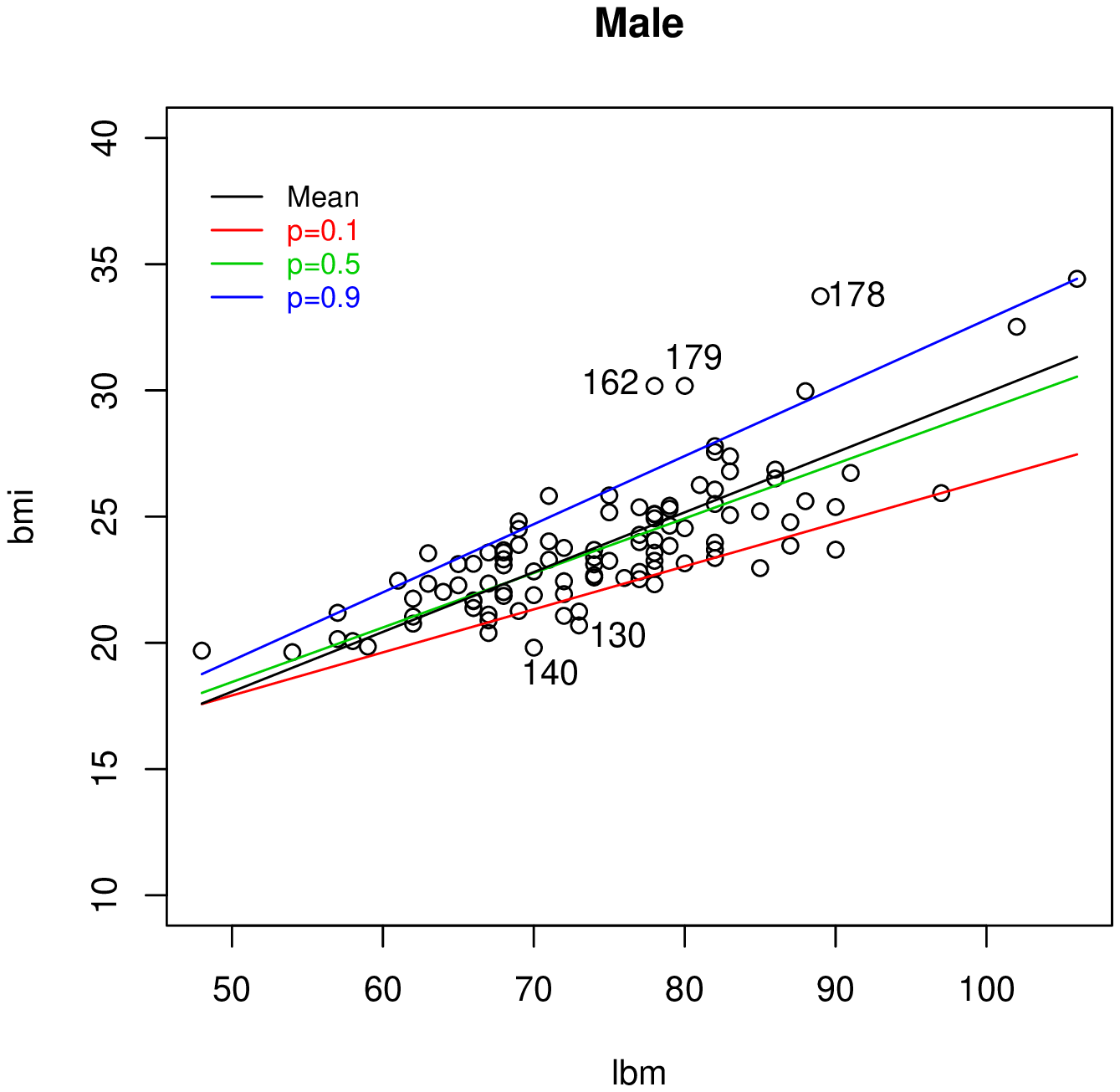}
\caption{AIS data: Fitted regression lines for the three selected
quantiles along with the mean regression line. The influential
observations are numbered. \label{fig:1b}}
\end{center}
\end{figure}

Figure \ref{fig:1b} shows the estimated quartiles of two levels of
gender at each LBM point from our EM algorithm along with the
estimates obtained via mean regression. From this figure we can see
clear attenuation in $\beta_1$ due to the use of the median
regression related to the mean regression. It is possible to observe
in this figure some atypical individuals that could have an
influence on the ML estimates for different values of quantiles.  In
this figure, the individuals $\#75, ~\#130, ~\#140 ~\#162, ~\#160$
and $\#178$ were marked since they were detected as potentially
influential.

\begin{figure}[h!]
\begin{center}
\centering \hspace {1cm}\centering \\
\includegraphics[scale=0.35]{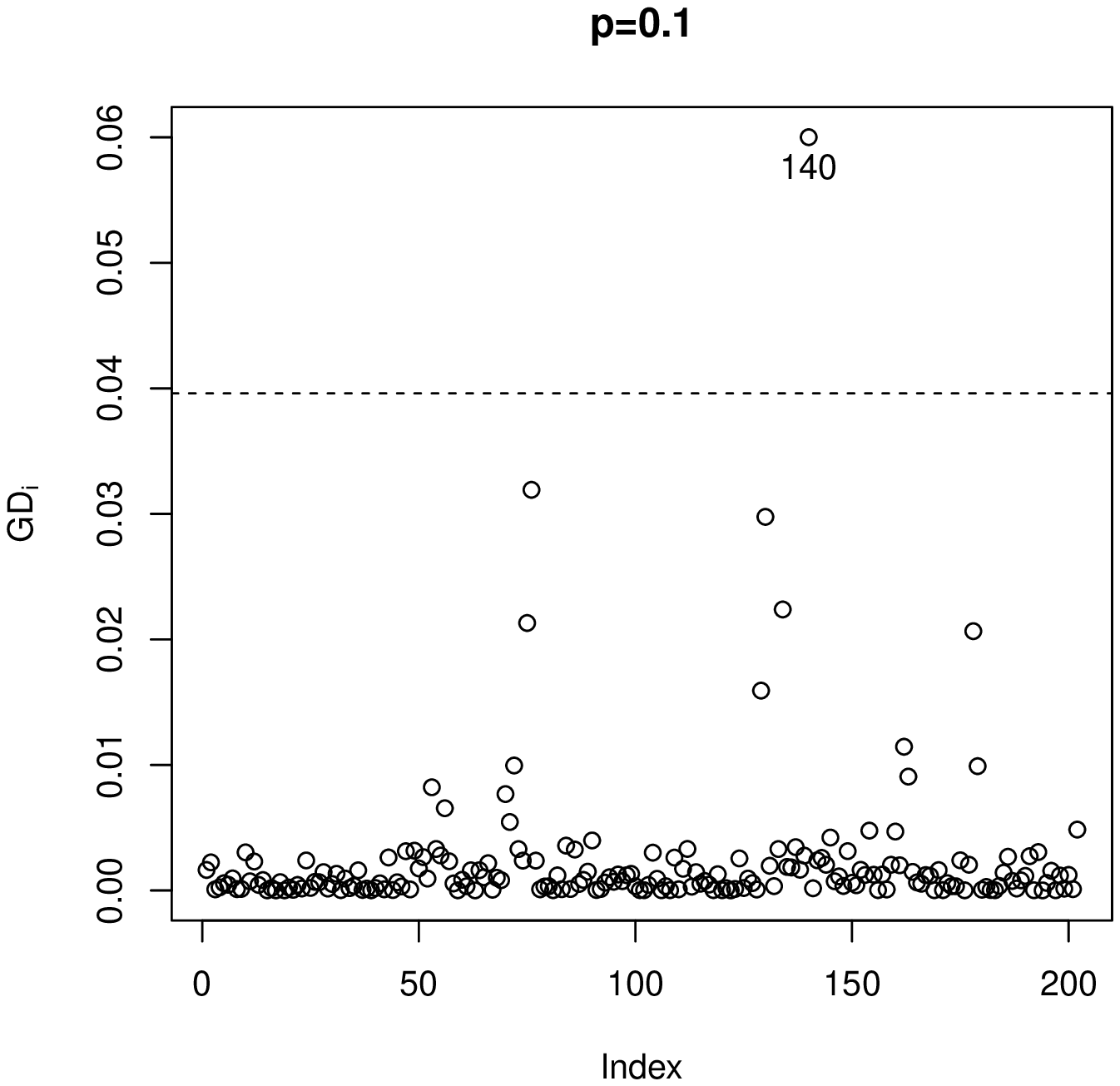}~\includegraphics[scale=0.35]{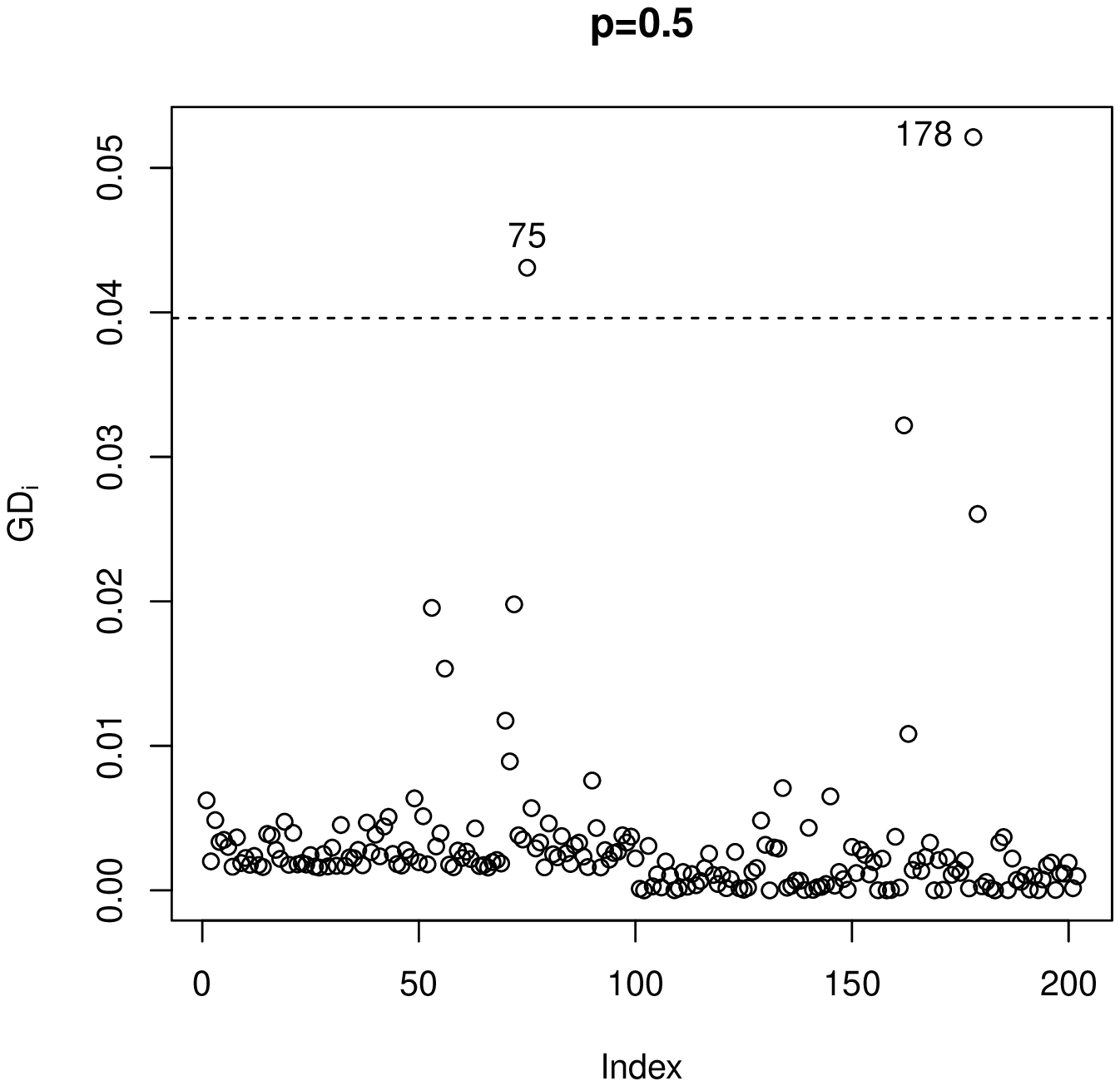}~\includegraphics[scale=0.35]{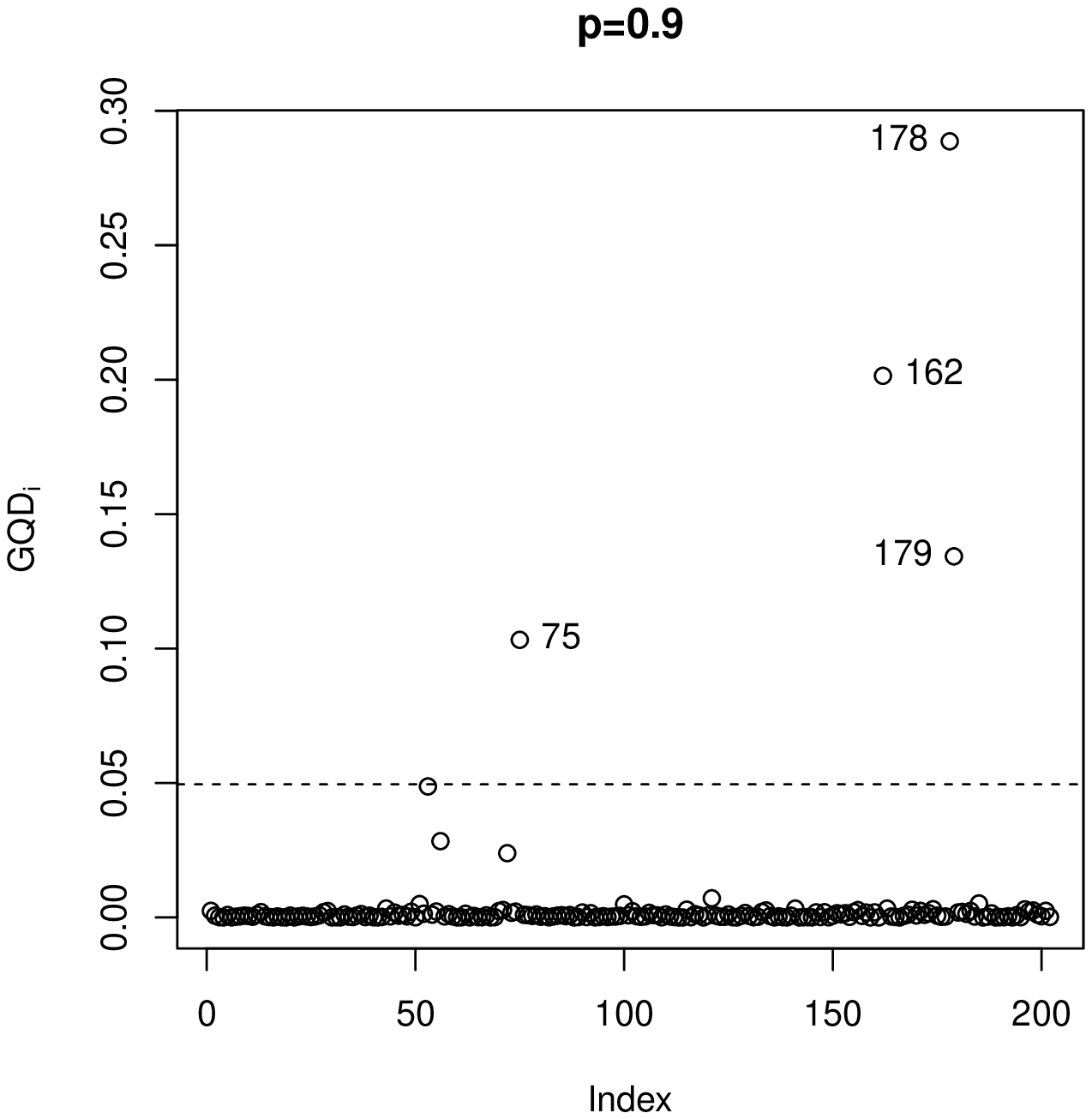}\\
\includegraphics[scale=0.35]{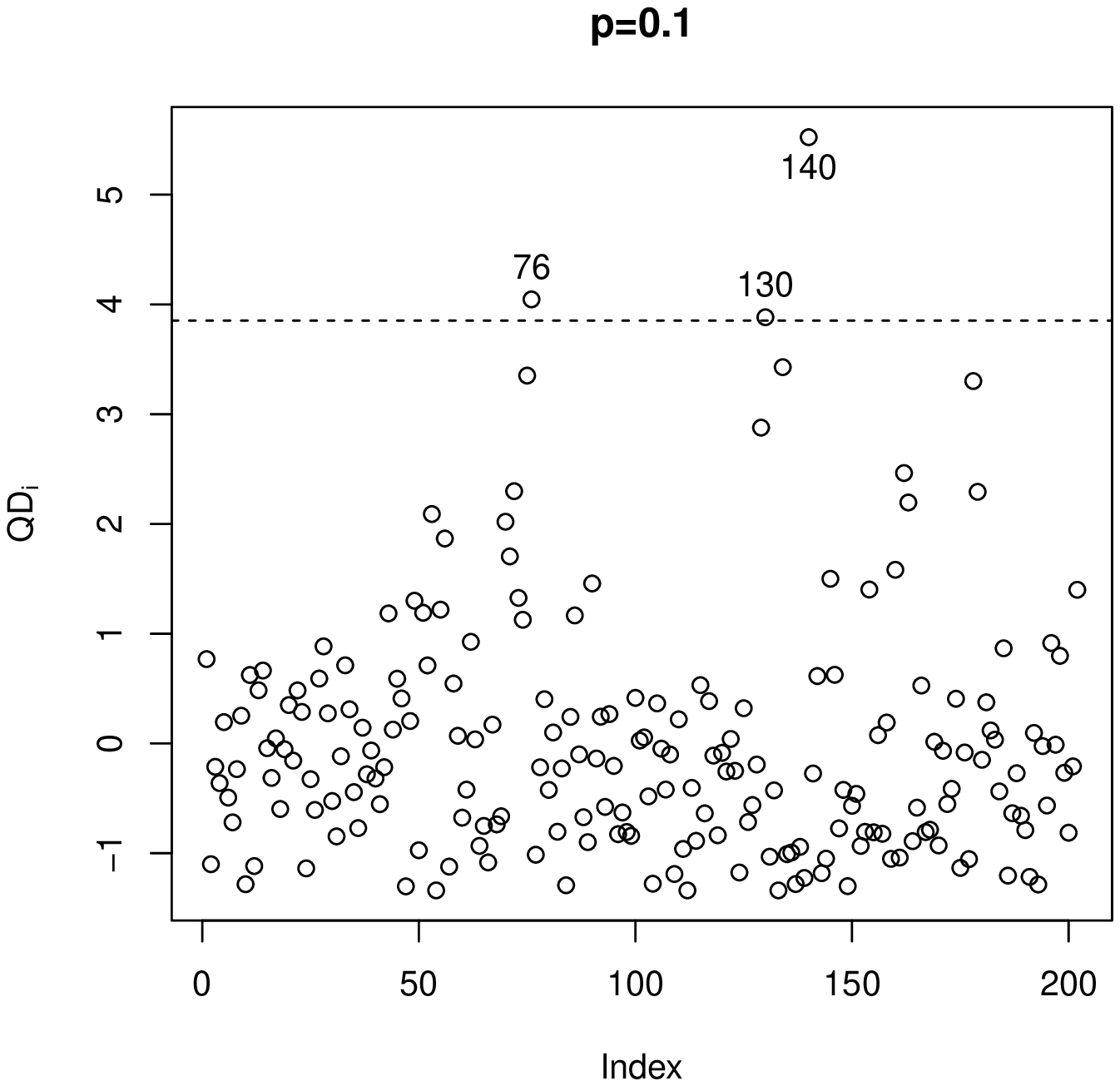}~\includegraphics[scale=0.35]{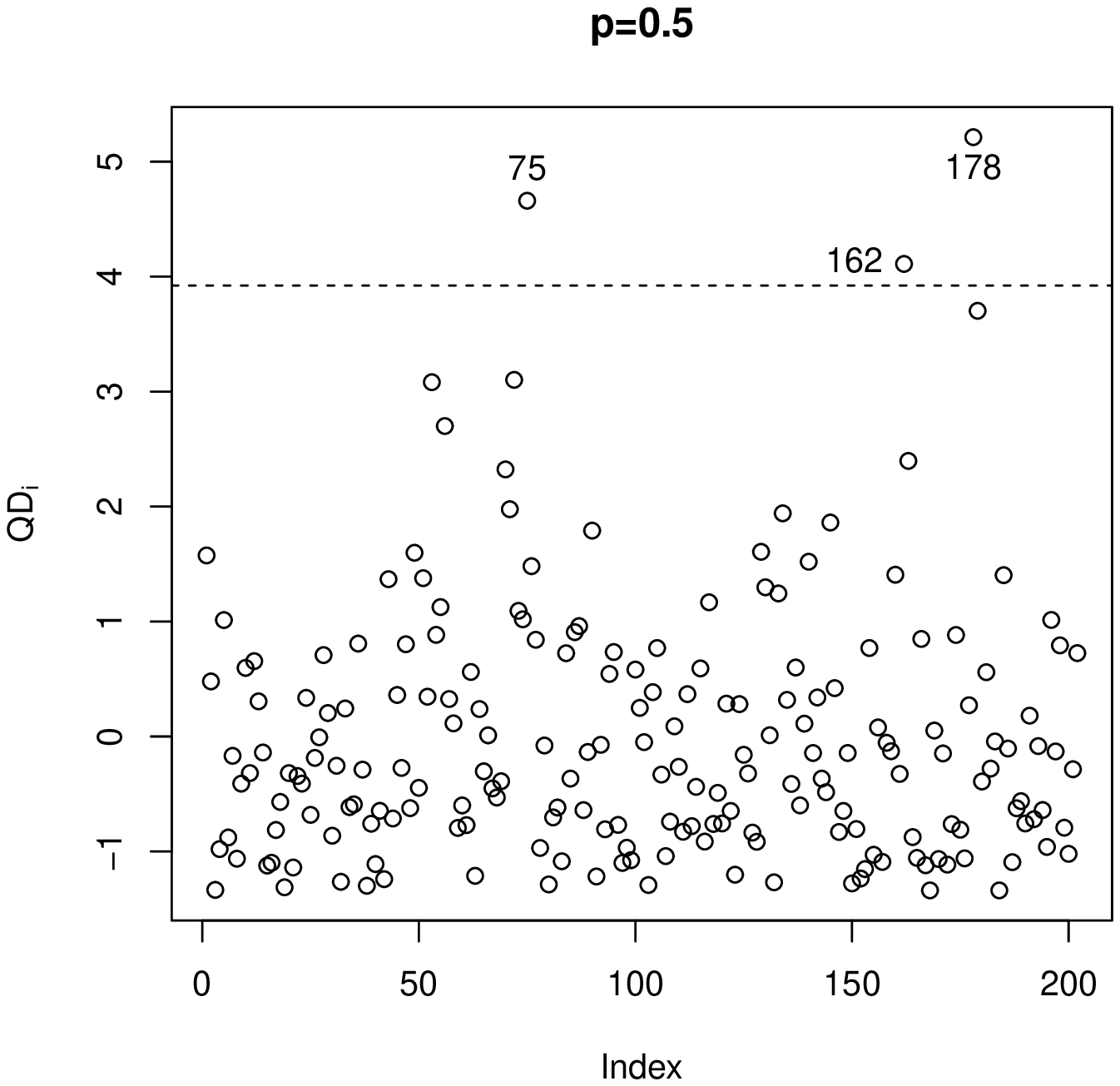}~\includegraphics[scale=0.35]{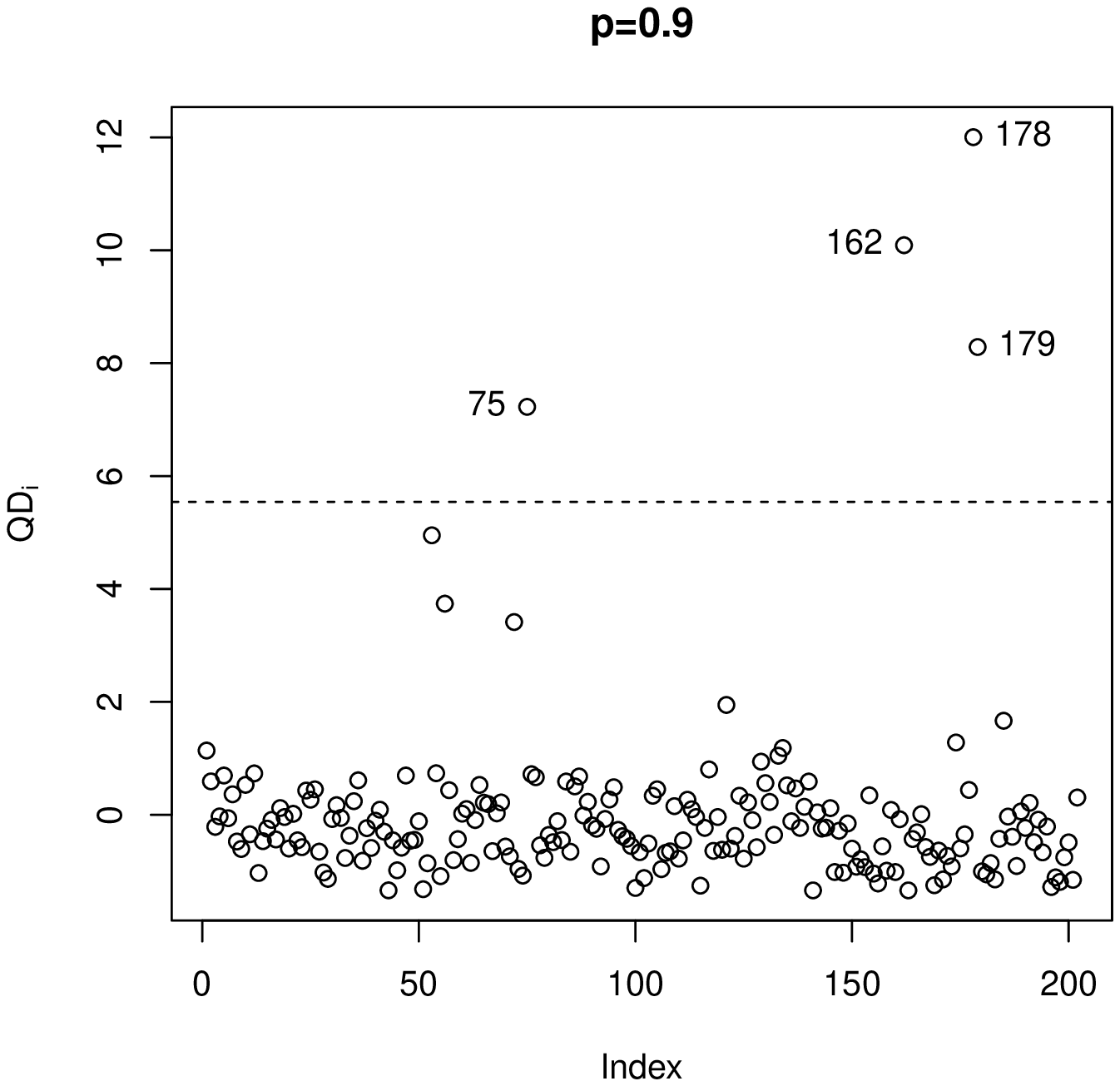}\\
\caption{Index plot of (first row) approximate likelihood distance
$GD^1_i$. (second row). Index plot of approximate likelihood
displacement $QD^1_i$. The influential observations are
numbered.\label{pert3} }
\end{center}
\end{figure}

In order to identify influential observations at different quantiles
when some observation is eliminated, we can generate graphs of the
generalized Cook distance $GD^{l}_i$, as explained in Section
\ref{Sec Diagnostic}. A high value for $GD^{l}_{i}$ indicates that
the $i$th observation has a high impact on the maximum likelihood
estimate of the parameters. Following \cite{Barros}, we can use
$2(p+1)/n$ as benchmark for the $GD^{l}_{i}$ at different quantiles.
Figure \ref{pert3} (first row) presents the index plots of
$GD^{l}_i$. We note from this figure that, only observation $\#140$
appears as influential in the ML estimates at $p=0.1$ and
observations $\#75, \#178$ as influential at $p=0.5$, whereas
observations $\#75,\#162, \#178$ and $\#179$ appear as influential
in the ML estimates at $p=0.9$. Figure \ref{pert3} (second row)
presents the index plots of $QD^1_i$. From this figure, it can be
noted that observations $\#76,\#130, \#140$ appear to be influential
at $p=0.1$, whereas observations $\#75,\#162$ and $\#178$ seem to be
 influential in the ML estimates at $p=0.1$, and in addition
observation  $\#179$ appears to be influential at $p=0.9$.

\section{Simulation studies} \label{sec simulation study}
In this section, the results from two simulation studies are
presented to illustrate the performance of the proposed method.

\subsection{Robustness of the EM estimates (Simulation study 1)}

We conducted a simulation study to assess the performance of the
proposed EM algorithm, by mimicking the setting of the AIS data by
taking the sample size $n = 202$. We simulated data from the model

\begin{equation}
y_{i}=\beta_{1} + \beta_{2}x_{i2} + \beta_{3}x_{i3} + \epsilon_{i},
\,\,\,\,\,\,\,\, i=1,\ldots,202, \label{simulation_1}
\end{equation}
where the $x_{ij}'s$ are simulated from a uniform distribution
(U(0,1)) and the errors $\epsilon_{ij}$ are simulated from four
different distributions: $(i)$ the standard normal distribution
$N(0,1)$, $(ii)$ a Student-t distribution with three degrees of
freedom, $t_{3}(0,1)$, $(iii)$ a heteroscedastic normal
distribution, $(1+x_{i2})N(0,1)$ and, $(iv)$ a bimodal mixture
distribution $0.6t_3(-20,1)+0.4t_3(15,1)$. The true values of the
regression parameters were taken as $\beta_1=\beta_2=\beta_3=1$. In
this way, we had four settings and for each setting we generated
$10000$ data sets.

Once the simulated data were generated, we fit a QR model, with $p=
0.1,\, 0.5$ and $0.9$, under Barrodale and Roberts (BR), Lasso
(Lasso) and EM algorithms by using the "quantreg()" package and our
\verb"ALDqr()" package, from the R language, respectively. For the
four scenarios, we computed the  bias and the square root of the
mean square error (RMSE), for each parameter over the $M=10,000$
replicas. They are defined as:
\begin{eqnarray}
Bias(\gamma) &=& \overline{\widehat{\gamma}}-\gamma\label{bias} \,\, {\rm and}\, \, \, \,
RMSE(\gamma) = \sqrt{SE(\gamma)^2 + Bias(\gamma)^2}\label{EQM}
\end{eqnarray}
where $\overline{\widehat{\gamma}} =
\frac{1}{M}\sum_{i=1}^{M}\widehat{\gamma}_i$ and  $SE(\gamma)^2 =
{\frac{1}{M-1}\sum_{i=1}^{M}\lp\widehat{\gamma}_i -
\overline{\widehat{\gamma}}\rp^2},$ with $\gamma =
\beta_1,\beta_2,\beta_3$ or $\sigma$, $\widehat{\gamma}_i$ is the
estimate of $\gamma$ obtained in replica $i$ and $\gamma$ is the
true value. Table \ref{table.simul1} reports the simulation results for $p =
0.1,\, 0.5$ and $0.9$. We observe that the EM yields lower biases
and RMSE than the other two estimation methods under all the
distributional scenarios. This finding suggests that the EM would
produce better results than other alternative methods typically used
in the literature of QR models.
\begin{table}[htbp!]
  \centering
  \caption{Simulation study. Bias and root mean-squared
  error (RMSE)  of $\bbeta$ under different error distributions. The estimates under Barrodale and Roberts (BR)
and Lasso  (Lasso) algorithms were obtained by the "quantreg()"
package from the R language.} {\footnotesize{
\begin{tabular}{lccccccc}

        \hline\hline
    \multicolumn{2}{c}{} & \multicolumn{2}{c}{$\beta_{1}$}  & \multicolumn{2}{c}{$\beta_{2}$} & \multicolumn{2}{c}{$\beta_{3}$} \\

     \hline
    Method                     & $p$ &    Bias &  RMSE  & Bias & RMSE & Bias & RMSE  \\
        \hline
    $\epsilon \sim N(0,1)$ & \multicolumn{7}{c}{} \\
        BR  & 0.1 & -1.2639 & 1.3444 &   0.0076 &  0.5961 &-0.0030 &  0.5934  \\
                               & 0.5 &  0.0064 & 0.3376 &  -0.0048 &  0.4390 &-0.0051 &  0.4453  \\
                               & 0.9 &  1.2640 & 1.3460 &   0.0030 &  0.6051 & 0.0069 &  0.6039 \\

               LPQR                  & 0.1 & -0.9664 & 1.0464 &  -0.3072 &  0.6165 &-0.3110 &  0.6187  \\
                               & 0.5 &  0.1474 & 0.3628 &  -0.1463 &  0.4534 &-0.1462 &  0.4576  \\
                               & 0.9 &  1.5901 & 1.6460 &  -0.3164 &  0.6173 &-0.3076 &  0.6179  \\

        EM                     & 0.1 & -1.2551 & 1.3362 &  -0.0055 &  0.5964 &-0.0090 &  0.6020  \\
                               & 0.5 &  0.0040 & 0.3286 &  -0.0050 &  0.4332 &-0.0031 &  0.4363  \\
                               & 0.9 &  1.2694 & 1.3484 &  -0.0071 &  0.6019 &-0.0120 &  0.5955  \\
                                                    \hline
    $\epsilon \sim t_{3}(0,1)$  & \multicolumn{7}{c}{} \\
        BR  & 0.1 & -1.2446 & 1.3364 &  -0.0290 &  0.6274 &-0.0313 &  0.6259 \\
                               & 0.5 &  0.1049 & 0.4870 &   0.1213 &  0.6714 & 0.1123 &  0.6708  \\
                               & 0.9 &  2.3618 & 2.8408 &   1.0056 &  2.4928 & 0.9459 &  2.4332 \\

                LPQR                  & 0.1 & -0.9315 & 1.0219 &  -0.3478 &  0.6422 &-0.3412 &  0.6354  \\
                               & 0.5 &  0.3007 & 0.5410 &  -0.0928 &  0.6310 &-0.0831 &  0.6237  \\
                               & 0.9 &  3.0443 & 3.2880 &   0.1911 &  1.6375 & 0.2231 &  1.6601 \\

        EM                     & 0.1 & -1.2287 & 1.3213 &  -0.0402 &  0.6209 &-0.0374 &  0.6265  \\
                               & 0.5 &  0.0965 & 0.4866 &   0.1352 &  0.6789 & 0.1304 &  0.6758  \\
                               & 0.9 &  2.3781 & 2.8459 &   0.9464 &  2.4082 & 0.9264 &  2.4167  \\
                                                    \hline
     $\epsilon \sim (1+x_{2})N(0,1)$  & \multicolumn{7}{c}{} \\
        BR  & 0.1 & -1.2869 & 1.4256 &   0.0130 &  0.8706 &-1.2554 &  1.5381  \\
                               & 0.5 & -0.0051 & 0.4468 &   0.0049 &  0.6336 & 0.0061 &  0.6509 \\
                               & 0.9 &  1.2868 & 1.4259 &   0.0018 &  0.8686 & 1.2307 &  1.5256 \\

                LPQR                  & 0.1 & -1.1393 & 1.2272 &  -0.3694 &  0.7773 &-1.1450 &  1.2756  \\
                               & 0.5 &  0.1834 & 0.4520 &  -0.1906 &  0.6193 &-0.1963 &  0.6304  \\
                               & 0.9 &  1.6972 & 1.7933 &  -0.3621 &  0.7925 & 0.7494 &  1.1587  \\

        EM                     & 0.1 & -1.2772 & 1.4140 &   0.0051 &  0.8646 &-1.2341 &  1.5195 \\
                               & 0.5 &  0.0954 & 0.4892 &   0.1289 &  0.6724 & 0.1316 &  0.6694 \\
                               & 0.9 &  1.2599 & 1.3987 &   0.0076 &  0.8723 & 1.2488 &  1.5315  \\
                                            \hline
       $\epsilon \sim 0.6t_3(-20,1)+0.4t_3(15,1)$  & \multicolumn{7}{c}{} \\
        BR  & 0.1 & -1.2350 & 1.3268 &  -0.0395 &  0.6160 &-0.0396 &  0.6192  \\
                               & 0.5 &  0.1029 & 0.4896 &   0.1214 &  0.6780 & 0.1212 &  0.6741 \\
                               & 0.9 &  2.3857 & 2.8737 &   0.9657 &  2.4574 & 0.9558 &  2.4585 \\

                LPQR                  & 0.1 & -0.9664 & 1.0464 &  -0.3072 &  0.6165 &-0.3110 &  0.6187  \\
                               & 0.5 &  0.1474 & 0.3628 &  -0.1463 &  0.4534 &-0.1462 &  0.4576  \\
                               & 0.9 &  1.5901 & 1.6460 &  -0.3164 &  0.6173 &-0.3076 &  0.6179  \\

        EM                     & 0.1 & -0.9327 & 1.0201 &  -0.3491 &  0.6433 &-0.3355 &  0.6372 \\
                               & 0.5 &  0.2880 & 0.5343 &  -0.0745 &  0.6216 &-0.0717 &  0.6159 \\
                               & 0.9 &  3.0624 & 3.3102 &   0.1702 &  1.6627 & 0.2221 &  1.6575  \\
                               \hline\hline

\end{tabular}
}}

  \label{table.simul1}%
\end{table}%

\subsection{Asymptotic properties (Simulation study 2)} \label{sec simulation study 2}

We also conducted a simulation study to evaluate the finite-sample
performance of the parameter estimates. We generated artificial
samples from the regression model (\ref{simulation_1}) with
$\beta_1=\beta_2=\beta_3=1$ and $x_{ij}\sim U(0,1)$. We chose
several distributions for the random term $\epsilon_i$ a little
different than the simulation study 1, say, $(i)$ normal
distribution $N(0,2)$ (N1), $(ii)$ a Student-t distribution
$t_{3}(0,2)$ (T1), $(iii)$ a heteroscedastic normal distribution,
$(1+x_{i2})N(0,2)$ (N2) and, $(iv)$ a bimodal mixture distribution
$0.6t_3(-20,2)+0.4t_3(15,2)$ (T2). Finally, the sample sizes were
fixed at $n = 50, 100, 150, 200, 300,$ $400, 500, 700$ and $800$.

\begin{figure}[!t]
\begin{center}
\includegraphics[scale=0.35]{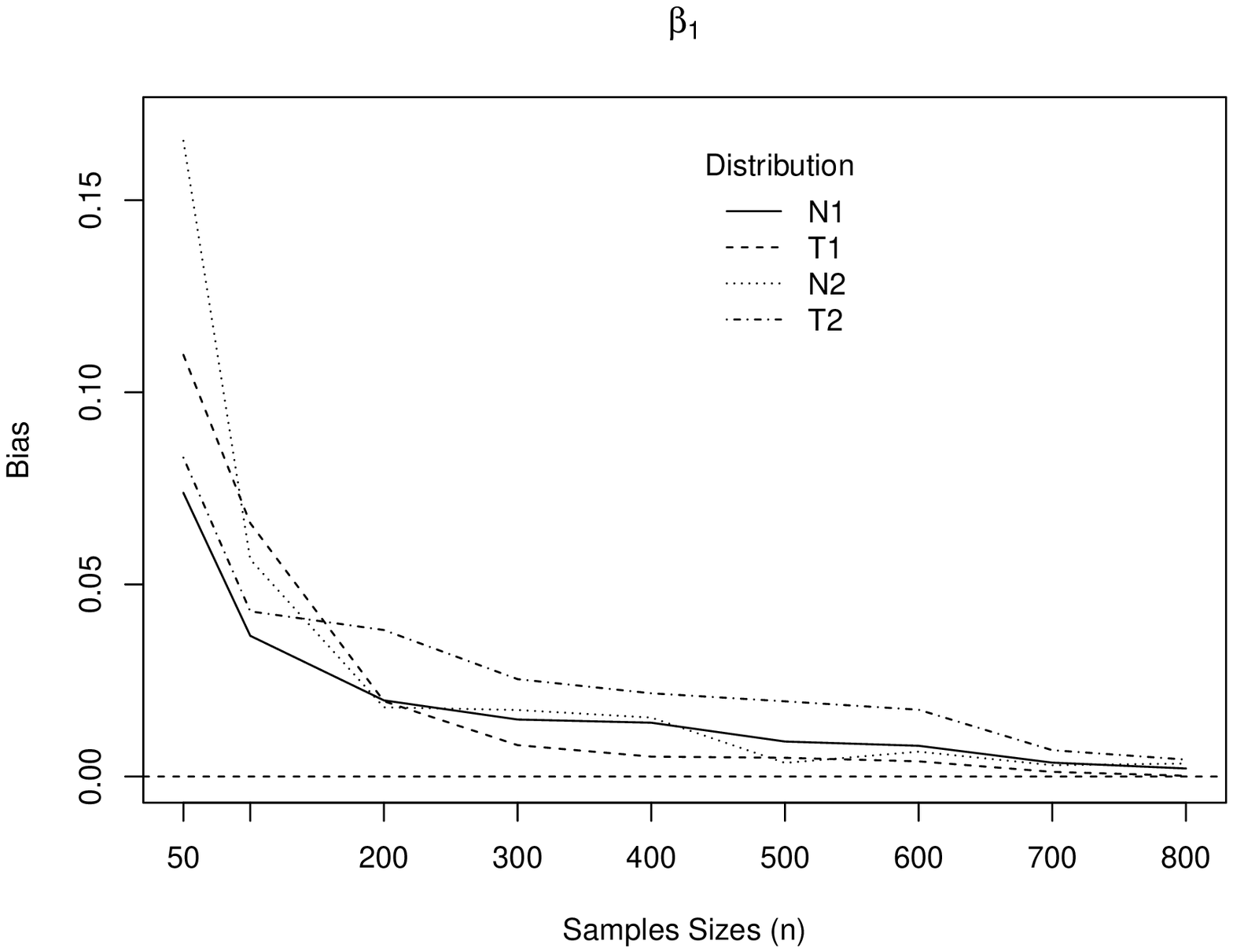}~\includegraphics[scale=0.35]{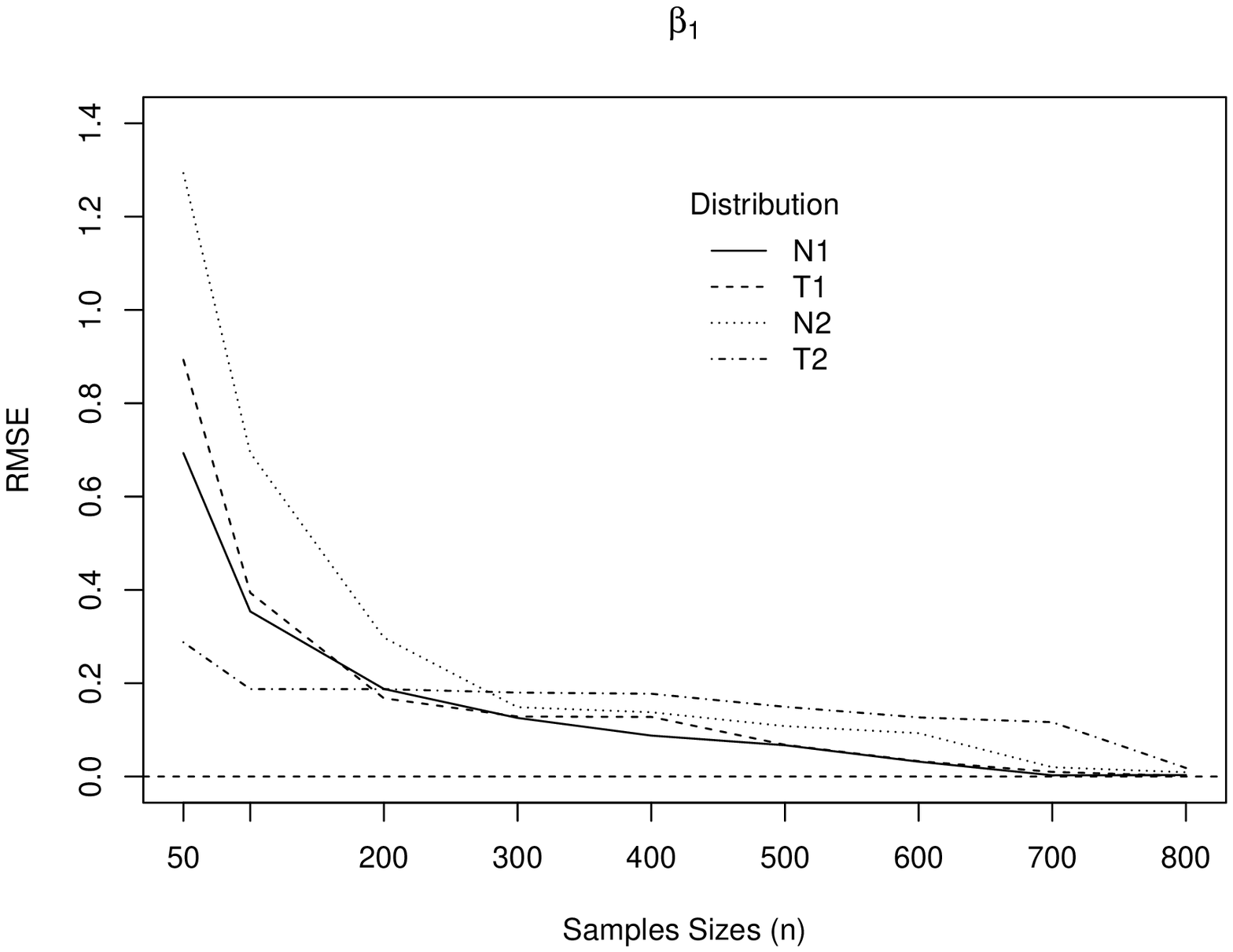}\\
\includegraphics[scale=0.35]{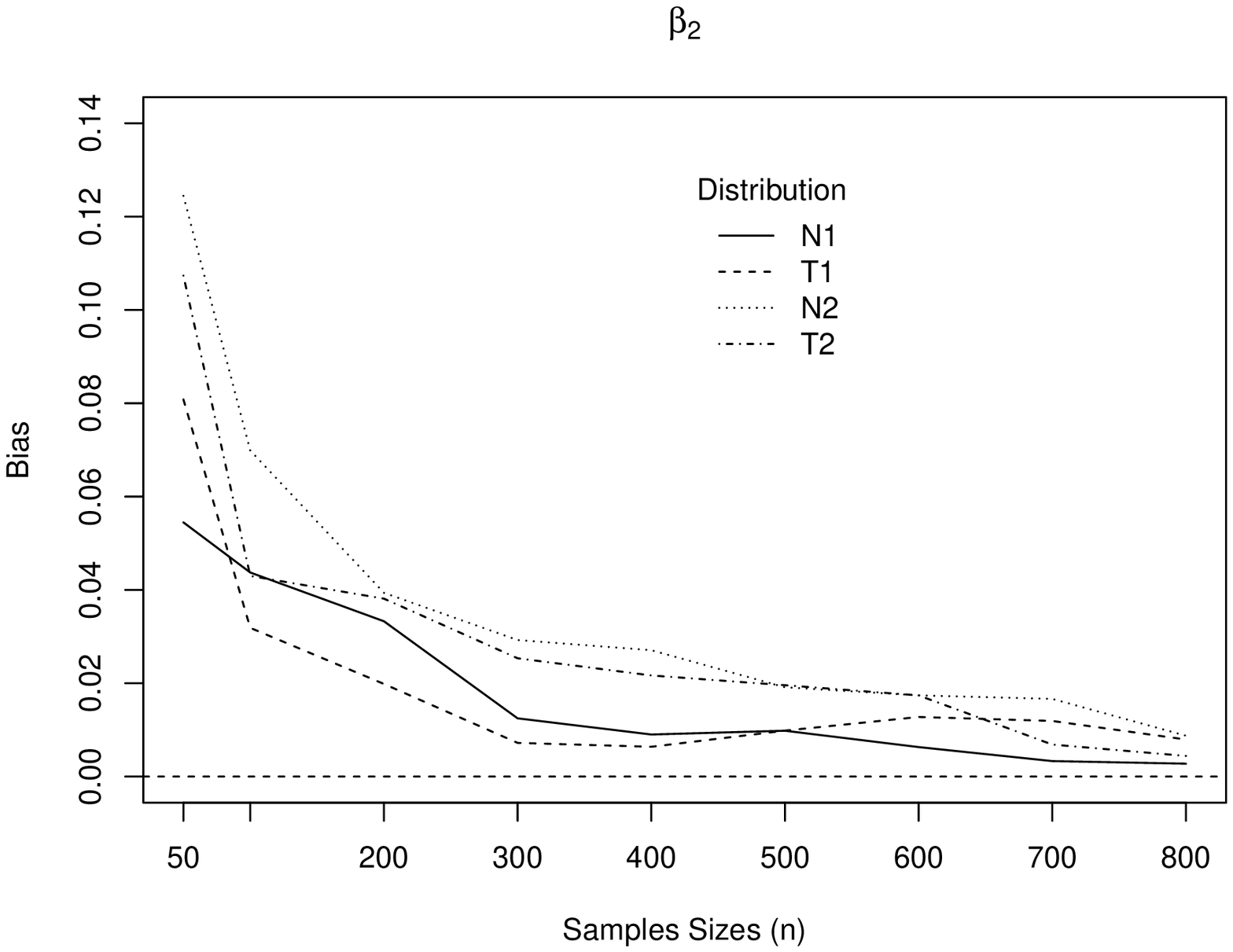}~\includegraphics[scale=0.35]{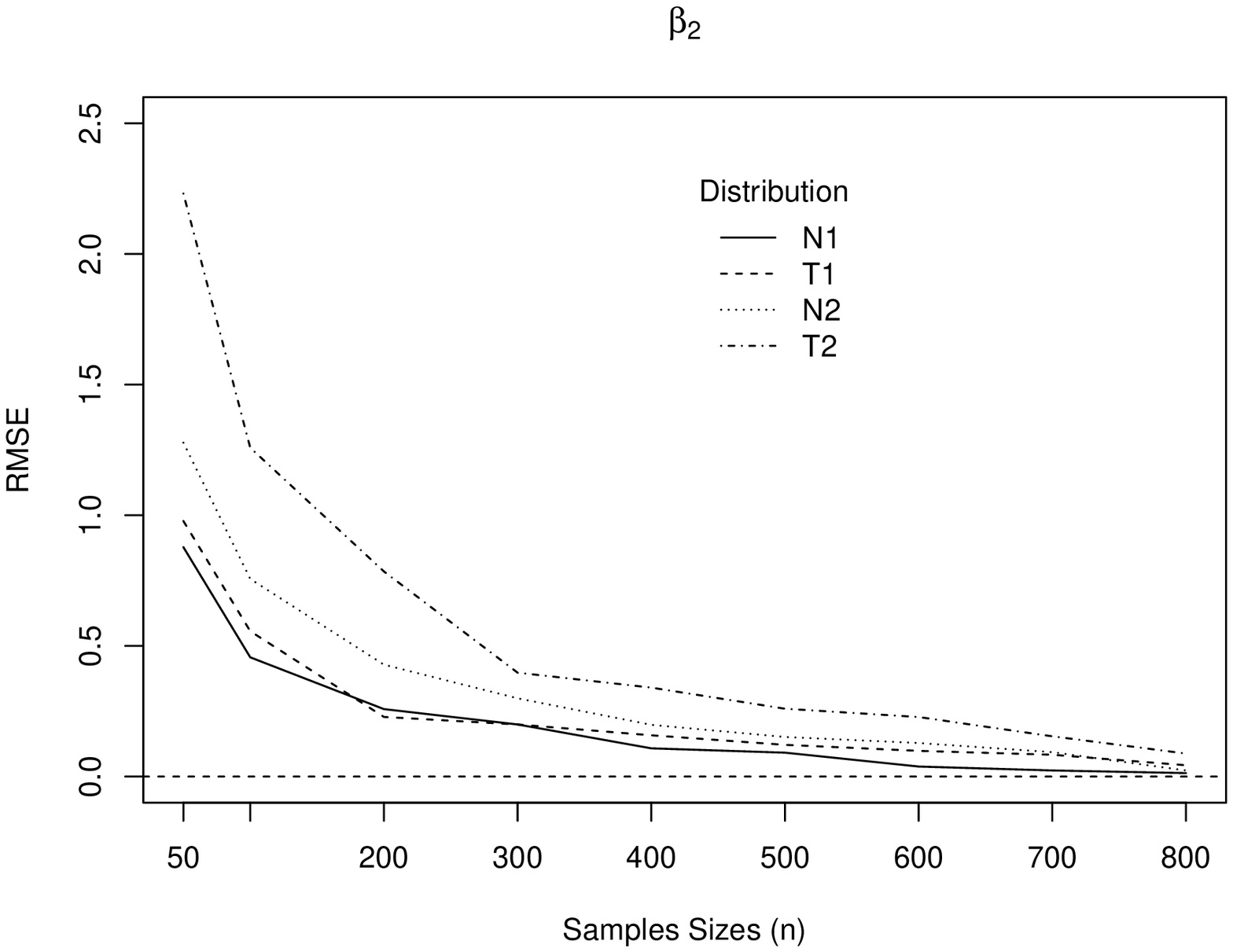}\\
\includegraphics[scale=0.35]{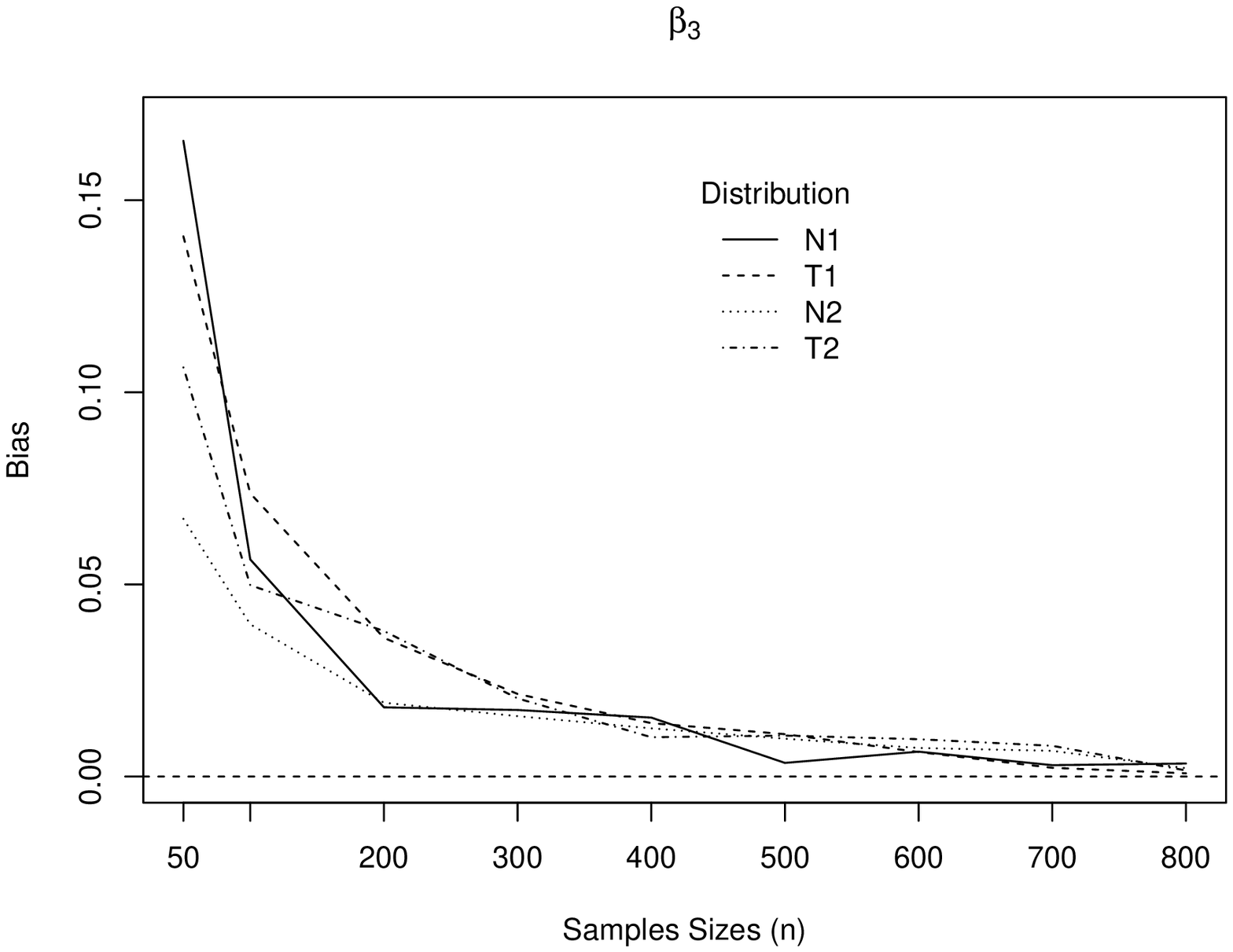}~\includegraphics[scale=0.35]{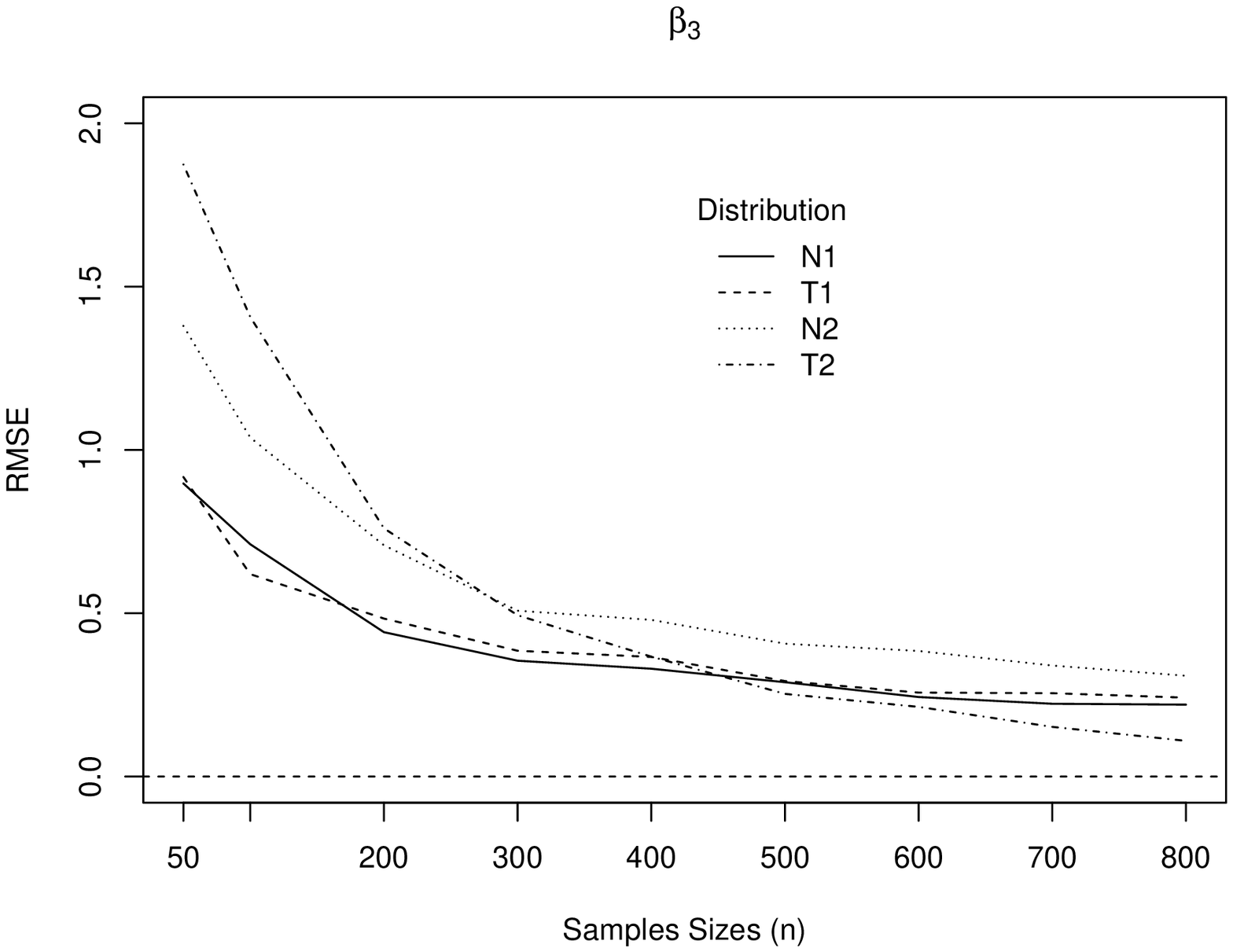}\\
\caption{Simulation study 2. Average bias (first column) and average
MSE (second column) of the estimates of $\beta_1$,$\beta_2$,
$\beta_3$ with $p=0.5$ (median regression), where $N1=N(0,2)$,
$T1=t_3(0,2)$, $N2=(1+x_2)N(0,2)$ and
$T2=0.6t_3(-20,2)+0.4t_3(15,2)$ .\label{fig:77a}}
\end{center}
\end{figure}

For each combination of parameters and sample sizes, $10000$ samples
were generated under the four different situations of error
distributions (N1, T1, N2, T2). Therefore, 36 different simulation
runs are performed. Once all the data were simulated, we fit the QR
model with $p=0.5$ and the  bias (\ref{bias}) and the square root of
the mean square error (\ref{EQM}) were recorded. The results  are
shown in Figure \ref{fig:77a}. We can see a pattern of convergence
to zero of the bias and MSE when $n$ increases. As a general rule,
we can say that bias and MSE tend to approach to zero when the
sample size increases, indicating that the estimates based on the
proposed EM-type algorithm do provide good asymptotic properties.
This same pattern of convergence to zero is repeated considering
different levels of the quantile $p$.

\section{Conclusion}
We have studied a likelihood-based approach to the estimation of the
QR based on the asymmetric Laplace distribution (ALD). By utilizing
the relationship between the QR check function and the ALD, we cast
the QR problem into the usual likelihood framework. The mixture
representation  of the ALD allows us to express a QR model as a
normal regression model, facilitating the implementation of an  EM
algorithm, which naturally provides the ML estimates of the model
parameters with the observed information matrix as a by product. The
EM algorithm  was implemented as part of the R package \textit{ALDqr()}. We hope that by making the code of our method
available, we will lower the barrier for other researchers to use
the EM algorithm  in their studies of quantile regression. Further,
we presented diagnostic analysis in QR models, which was based on
the case-deletion technique suggested by \cite{zhu2001case} and
\cite{ZhuLee2001},  which are the counterparts for missing data
models of the well-known ones proposed by \cite{cook77} and
\cite{cook86}. The simulation studies  demonstrated the superiority
of the proposed methods to the existing methods, implemented in the
package \verb"quantreg()". We applied our methods to a real data set
(freely downloadable from R) in order to illustrate how the
procedures can be used to identify outliers and to obtain robust ML
parameter estimates. From these results, it is encouraging that the
use of ALD offers a better alternative in the analysis of QR models.

Finally, the  proposed  methods can  be extended to a more general
framework, such as, censored (Tobit) regression models, measurement
error models, nonlinear regression models, stochastic volatility
models, etc and should yield satisfactory results at the expense of
additional complexity in implementation. An in-depth investigation
of such extensions is beyond the scope of the present paper, but
these are interesting topics for further research.

\section*{Acknowledgements}
The research of V. H. Lachos was supported by Grant 305054/2011-2
from Conselho Nacional de Desenvolvimento Cient\'{\i}fico e Tecnol\'{o}gico
(CNPq-Brazil) and by Grant 2014/02938-9 from  Funda\c c\~ao de
Amparo \`{a} Pesquisa do Estado de S\~ao Paulo (FAPESP-Brazil).



\end{document}